%
%
%
%
\magnification=1200
\input amstex
\input amsppt.sty
\def\SBIMSMark#1#2#3{
 \font\SBF=cmss10 at 10 true pt
 \font\SBI=cmssi10 at 10 true pt
 \setbox0=\hbox{\SBF Stony Brook IMS Preprint \##1}
 \setbox2=\hbox to \wd0{\hfil \SBI #2}
 \setbox4=\hbox to \wd0{\hfil \SBI #3}
 \setbox6=\hbox to \wd0{\hss
             \vbox{\hsize=\wd0 \parskip=0pt \baselineskip=10 true pt
                   \copy0 \break%
                   \copy2 \break%
                   \copy4 \break}}
 \dimen0=\ht6   \advance\dimen0 by \vsize \advance\dimen0 by 8 true pt
                \advance\dimen0 by -\pagetotal
 \dimen2=\hsize \advance\dimen2 by .25 true in
%
%
  \openin2=publishd.tex
  \ifeof2\setbox0=\hbox to 0pt{}
  \else 
     \setbox0=\hbox to 3.1 true in{
                \vbox to \ht6{\hsize=3 true in \parskip=0pt  \noindent  
                \input publishd.tex 
                \vfill}}
  \fi
  \closein2
  \ht0=0pt \dp0=0pt
 \ht6=0pt \dp6=0pt
 \setbox8=\vbox to \dimen0{\vfill \hbox to \dimen2{\copy0 \hss \copy6}}
 \ht8=0pt \dp8=0pt \wd8=0pt
 \copy8
 \message{*** Stony Brook IMS Preprint #1, #2 ***}
}

\input psbox
\psfordvips
\let\fillinggrid=\relax
\UseAMSsymbols
\documentstyle{amsppt}
\refstyle{C}
\NoBlackBoxes
\NoRunningHeads
\vsize=21.2truecm
\hsize=15 truecm
\hcorrection{0.8truecm}
\def\today{November, 1997}

\newcount\refno\refno=1
\define\referencelabel{\number\refno} 
\define\newrefno{\advance\refno by1}
\edef\refnoAh
{\referencelabel}\newrefno
\define\refAh{{\bf[\refnoAh]}}
\edef\refnodF
{\referencelabel}\newrefno
\define\refdF{{\bf[\refnodF]}}
\edef\refnodFdM
{\referencelabel}\newrefno
\define\refdFdM{{\bf[\refnodFdM]}}
\edef\refnoGrS
{\referencelabel}\newrefno
\define\refGrS{{\bf[\refnoGrS]}}
\edef\refnoHea
{\referencelabel}\newrefno
\define\refHea{{\bf[\refnoHea]}}
\edef\refnoHeb
{\referencelabel}\newrefno
\define\refHeb{{\bf[\refnoHeb]}}
\edef\refnoLYam
{\referencelabel}\newrefno
\define\refLYam{{\bf[\refnoLYam]}}
\edef\refnoMca
{\referencelabel}\newrefno
\define\refMca{{\bf[\refnoMca]}}
\edef\refnoMcb
{\referencelabel}\newrefno
\define\refMcb{{\bf[\refnoMcb]}}
\edef\refnoMcc
{\referencelabel}\newrefno
\define\refMcc{{\bf[\refnoMcc]}}
\edef\refnoMS
{\referencelabel}\newrefno
\define\refMS{{\bf[\refnoMS]}}
\edef\refnoSu
{\referencelabel}\newrefno
\define\refSu{{\bf[\refnoSu]}}
\edef\refnoSwa
{\referencelabel}\newrefno
\define\refSwa{{\bf[\refnoSwa]}}
\edef\refnoYam
{\referencelabel}\newrefno
\define\refYam{{\bf[\refnoYam]}}
\edef\refnoYo
{\referencelabel}\newrefno
\define\refYo{{\bf[\refnoYo]}}
%
%
\define\AcM{Acta Math.}

\define\CMP{Commun. Math. Phys.}
\define\ETDS{Ergod. Th. \& Dynam. Sys.}
\define\IHES{Publ. Math. IHES}

%
%
\define\Ah{L.~Ahlfors}

\define\dF{E.~de~Faria}

\define\Gk{J.~Graczyk}
\define\He{M.~Herman}

\define\Ly{M.~Lyubich}
\define\Mc{C.~McMullen}
\define\Me{W.~de~Melo}

\define\Su{D.~Sullivan}
\define\Sw{G.~\'Swi\c{a}tek}

\define\vS{S.~van~Strien}
\define\Yam{M.~Yampolsky}
\define\Yo{J.~C.~Yoccoz}
\overfullrule=0pt

\def\diam{\roman{diam\,}}
\def\dist{\roman{dist\,}}
\def\dom{\roman{Dom\,}}
\def\mod{\roman{mod\,}}
\def\Re{\roman{Re\,}}
\def\Im{\roman{Im\,}}
\define\({\left(}
\define\){\right)}  
\define\[{\left[}
\define\]{\right]}

\let\epsilon=\varepsilon
\newcount\eqnocounter
\eqnocounter=0
\define\slabel#1{\global\advance\eqnocounter by1%
\global\expandafter\edef\csname eqno#1\endcsname{(\number\eqnocounter)}%
\eqno({\number\eqnocounter})}
\define\alabel#1{\global\advance\eqnocounter by1%
\global\expandafter\edef\csname eqno#1\endcsname{(\number\eqnocounter)}%
({\number\eqnocounter})}
\define\endclaim{\egroup\par\ifdim\lastskip<\smallskipamount
\removelastskip\penalty20\smallskip\fi}
\define\claim{\smallbreak\noindent  
{\it Claim.\enspace}\ignorespaces\bgroup\sl}
\define\endassertion{\egroup\par\ifdim\lastskip<\smallskipamount
\removelastskip\penalty20\smallskip\fi}
\define\assertion(#1){\smallbreak\noindent 
{\rm(}{#1}{\rm)}\enspace\ignorespaces\bgroup\sl}
\define\endcase{\par\ifdim\lastskip<\smallskipamount
\removelastskip\penalty3000\smallskip\fi}
\define\case#1.{\smallbreak\noindent 
{\sl Case\/} #1.\enspace\ignorespaces}
\define\endcasebody{\par\ifdim\lastskip<\smallskipamount
\removelastskip\penalty-20\smallskip\fi}

%

\newcount\statemno\statemno=1
\define\newstatemno{\advance\statemno by1}
%
%
%
%
%
%

%
%
%
%
%
%
%
%

\nologo
\SBIMSMark{1997/17}{November 1997}{}

\topmatter

\title Rigidity of critical circle mappings II
\endtitle

\author Edson de Faria and Welington de Melo
\endauthor

\abstract We prove that any two real-analytic critical circle maps
with cubic critical point and the same irrational rotation number of bounded
type are $C^{1+\alpha}$ conjugate for some $\alpha>0$. 
\endabstract

\date
\today
\enddate

\address
Instituto de Matem\'atica e Estat\'\i stica,
Universidade de S\~ao Paulo,\endgraf
Rua do Mat\~ao 1010, CEP05508-900 S\~ao Paulo SP - Brasil
\endaddress 
\email edson\@ime.usp.br\endemail
~
\address
Instituto de Matem\'atica Pura e Aplicada, 
Estrada Dona Castorina 110, \endgraf
Jardim Bot\^anico, CEP22460-320
Rio de Janeiro RJ - Brasil
\endaddress
\email demelo\@impa.br\endemail

\keywords Holomorphic pairs, complex bounds, uniform twist, rigidity
\endkeywords 

\thanks This work has been partially supported by the Pronex Project on
Dynamical Systems,\endgraf by FAPESP Grant 95/3187-4 and by CNPq Grant
30.1244/86-3. 
\endthanks

\subjclass
Primary 58F03; Secondary 30F60, 58F23, 32G05
\endsubjclass

\endtopmatter


\heading
1. Introduction
\endheading

A celebrated rigidity theorem of M.~Herman {\refHea} states that any
two smooth diffeomorphisms of the circle with the same Diophantine
rotation number are smoothly conjugate.
Herman used purely real-variable techniques to prove this theorem.
In recent years, other fundamental rigidity results have been
established in the context of one-dimensional dynamics with the help
of complex-analytic techniques. Thus, in the context of circle
diffeomorphisms, J-C.~Yoccoz~{\refYo} used conformal surgery to describe the
optimal class of rotation numbers for which Herman's theorem holds in the
real-analytic category.  
In the context of interval maps, D.~Sullivan~{\refSu}
and C.~McMullen~{\refMcb} used a combination of real-variable
techniques with deep results from Teichm\"uller theory 
and hyperbolic geometry to establish the rigidity of the Cantor
attractor of a unimodal map with bounded combinatorics. The conjugacy
between two such maps with the same combinatorics is shown to be
$C^{1+\alpha}$ for some $\alpha>0$ depending only on the
combinatorics. In sharp contrast with the case of circle
diffeomorphisms, in this case one does not expect the conjugacy to be
much better than $C^1$, even when the unimodal maps are analytic.

In this paper we prove the following strong rigidity theorem for
real-analytic critical circle maps. 

\proclaim{Theorem {1.1}} Any two real-analytic critical circle maps with
the same irrational rotation number of bounded type are
$C^{1+\alpha}$ conjugate for some $0<\alpha<1$ depending only on the least
upper-bound of the coefficients of the continued fraction expansion of
$\rho$.
\endproclaim

The proof of this theorem uses both real and complex-analytic tools.
We use the {\it real a-priori bounds} of \'Swi\c{a}tek and Herman (see
{\refSwa}, {\refHeb} and {\refdFdM}) in the form stated in \S 2. The
complex-analytic 
part we develop here combines some of the powerful new ideas on 
conformal rigidity and renormalization recently developed by McMullen
in {\refMca}, {\refMcb}, {\refMcc}, with the basic theory of
holomorphic commuting pairs introduced in {\refdF} and a
generalization of the Lyubich-Yampolsky approach to the complex bounds
given in {\refLYam}, {\refYam}.

The general outline of the paper is as follows.
First we prove that after a {\it finite} number of renormalizations,
every real-analytic critical circle map with arbitrary irrational rotation
number gives rise to a {\it holomorphic commuting
pair} with good geometric bounds. 
Next, we prove that such a holomorphic pair (with good geometric
control) possesses McMullen's {\it uniform twisting} property
{\refMcb}. We also prove that its critical point is {\it
$\delta$-deep} for some universal $\delta >0$.
The general machinery developed by McMullen in {\refMca}, {\refMcb}
permits us to deduce from these facts that the conjugacy between any
two real-analytic critical circle maps with rotation number of type bounded by
$N$ is $C^{1+\varepsilon(N)}$-conformal at the critical point.
This implies that the successive renormalizations of both circle maps
around their critical points are converging together {\it
exponentially fast} in the $C^0$ sense (even $C^{1+\beta}$ for some
$\beta>0$), with rate of convergence given by a universal function of
$\varepsilon$. By a general proposition proved in {\refdFdM}, this implies
at last that the conjugacy between both maps is $C^{1+\alpha(N)}$
for some $\alpha(N)>0$. Although in our proof the H\"older exponent depends on
$N$, we believe that this should not be the case. 

It is to be noted that in {\refdFdM} we constructed $C^\infty$ counterexamples
(even with Diophantine rotation numbers) to this form of
$C^{1+\alpha}$ rigidity. One may ask whether such counterexamples exist in the
analytic category. We do not know the answer to this question. Nevertheless,
we conjecture that in the general case the conjugacy is always $C^{1+\alpha}$
at the critical point, for some universal $\alpha$ independent of the rotation
number. A possible approach to this conjecture is to use M.~Lyubich's recent
work on quadratic-like maps and universality of unimodal maps, but it is still
unclear at this time how to adapt his methods to holomorphic pairs and
critical circle maps. 

\subhead
Acknowledgments
\endsubhead
We wish to thank D.~Sullivan for sharing with us some of his deep
insights on renormalization.
We are grateful to C.~McMullen for explaining to us various aspects of 
his beautiful results on renormalization and rigidity.

\heading 
2. Preliminaries
\endheading
In this section, we introduce some basic concepts and notations.
We denote by $|I|$ the length of an interval $I$ on the line or the
circle. We use $\dist(\cdot,\cdot)$ for Euclidean distance, and
$\diam(\cdot)$ for Euclidean diameter.
All bounds achieved in this paper will depend on certain {\it
a-priori} constants (which can ultimately be traced back to the {\it
real a-priori bounds} given by Theorem {2.1} below). In this context, 
two positive quantities $a$ and $b$ are said to be {\it comparable} if
there exists a constant $C>1$ depending only on the {a-priori} bounds
such that $C^{-1}b\le a\le Cb$. We write $a\asymp b$ to denote that $a$
is comparable to $b$.

\subhead
Critical circle maps
\endsubhead
We identify the unit circle $S^1$ with the one-dimensional torus
${\Bbb R}/{\Bbb Z}$. A homeomorphism $f:S^1\to S^1$ which is at least
$C^1$ and satisfies $f'(c)=0$ and $f'(x)\neq 0$ for all $x\neq c$
is called a {\it critical circle map}. The point $c$ is the {\it
critical point} of $f$. In this paper $f$ will be real-analytic, and
the critical point will be cubic: this means that near $c$ we have
$f=\phi\circ Q\circ \psi$, where $\phi$, $\psi$ are real analytic
diffeomorphisms and $Q$ is the map $x\mapsto x^3$. 

We are interested in the geometry of orbits of $f$ only when there are 
no periodic points. In this case the rotation number of $f$ is
irrational and can be represented as an infinite continued fraction  
$$
\rho(f)\;=\;[a_0,a_1,\ldots,a_n,\ldots ]\;=\;
\cfrac 1\\ 
a_0+\cfrac 1\\ 
a_1+\cfrac 1\\ 
{}\cfrac {\cdots}\\
a_n+\cfrac 1\\ 
{\cdots}
\endcfrac
\ .
$$
When the partial quotients $a_n$ are bounded, we say that $\rho(f)$ is
a number of {\it bounded type}. We also refer to $\sup{a_n}$ as the
{\it combinatorial type} of $\rho(f)$. 

The denominators of the convergents of $\rho(f)$, defined recursively
by $q_0=1$, $q_1=a_0$ and $q_{n+1}=a_nq_n+q_{n-1}$ for all $n\ge 1$,
are the {\it closest return times} of the orbit of any point to
itself. 
We denote by $\Delta_n$ the closed interval
containing $c$ whose endpoints are $f^{q_n}(c)$ and $f^{q_{n+1}}(c)$.
We also let $I_n\subseteq \Delta_n$ be the closed interval whose
endpoints are $c$ and $f^{q_n}(c)$. Observe that
$\Delta_n=I_n\cup I_{n+1}$. The most important fact to remember when
studying the geometry of a circle map is that for each $n$ the
collection of intervals
$$
{\Cal P}_n\;=\;\Big\{I_n,\,f(I_n),\ldots,\, f^{q_{n+1}-1}(I_n)\Big\}
\cup\Big\{I_{n+1},\,f(I_{n+1}),\ldots ,\,f^{q_n-1}(I_{n+1})\Big\} 
$$
constitutes a partition of the circle (modulo endpoints), called {\it
dynamical partition of level $n$} of the map $f$. Note that, for all
$n$, ${\Cal P}_{n+1}$ is a refinement of ${\Cal P}_n$.

Of course, these definitions make sense for an arbitrary homeomorphism
of the circle. For a rigid rotation, we have
$|I_n|=a_{n+1}|I_{n+1}|+|I_{n+2}|$. Therefore, if $a_{n+1}$ is very
large then $I_n$ is much longer than $I_{n+1}$. It is a remarkable
fact, first proved by \'Swi\c{a}tek and Herman, that this never
happens for a critical circle map! The dynamical partitions ${\Cal
P}_n$ have {\it bounded geometry}, in the sense that adjacent
atoms have comparable lengths.

\proclaim{Theorem {2.1}} (The real bounds) Let $f:S^1\to S^1$ be a
critical circle map with arbitrary irrational rotation number. There
exists $n_0=n_0(f)$ such that for all $n\ge n_0$ and every pair $I,J$
of adjacent atoms of ${\Cal P}_n$ we have $K^{-1}|J|\le |I|\le K|J|$,
where $K>1$ is a universal constant.
\endproclaim

These {\it real a-priori bounds} are essential in all the
estimates that we perform in this paper. For proofs, see {\refGrS} and \S~3
of {\refdFdM}.
An important consequence of this theorem is the fact that every
critical circle map with rotation number of bounded type is
conjugate to the rotation with the same rotation number by a {\it
quasisymmetric} homeomorphism. Another consequence is the following
lemma.

\proclaim{Lemma {2.2}} Let $f$ be a critical circle map with arbitrary
rotation number, let $n\ge 1$, and let $J_{-i}=f^{q_{n+1}-i}(I_n)$ for $0\le
i<q_{n+1}$. Given $m<n$, let $i_1<i_2<\cdots<i_\ell$ be the moments in the
backward orbit $\{J_{-i}\}$ before the the first return to $I_{m+1}$ such that
$J_{-i_k}\subseteq I_m$. Then $\ell=a_{m+1}$, and we have 
$$
J_{-i_k}\subseteq f^{q_m+({a_{m+1}-k+1})q_{m+1}}(I_{m+1})
\ .
$$
Moreover, given an integer $M\ge 1$, there exists $C_M>1$
such that for all sufficiently large $n$ we have $C_M^{-1}|I_n|\le
|J_{-i_k}|\le C_M |I_n|$, provided $1\le k\le M$ or $a_{m+1}-M+1\le k\le
a_{m+1}$.
\endproclaim
\demo{Proof} The largest $j<q_{n+1}$ such that $f^j(I_n)\subseteq I_{m+1}$ is
easily computed as $j=q_{n+1}-q_{m+2}$. Since $q_{m+2}=q_m+a_{m+1}q_{m+1}$,
there are exactly $a_{m+1}$ subsequent moments $j<i<q_{n+1}$ such that
$f^i(I_n)\subseteq I_m$. The rest follows from Theorem {2.1} and the Koebe
distortion principle (see {\refMS}).\qed
\enddemo

\subhead
Commuting pairs and renormalization
\endsubhead
Let $f$ be a critical circle map as before, and let $n\ge 1$. The
first return map $f_n:\Delta_n\to \Delta_n$ to $\Delta_n=I_n\cup
I_{n+1}$, called the {\it $n$-th renormalization of $f$ without
rescaling}, is determined by a pair of maps, namely
$\xi=f^{q_n}:I_{n+1}\to \Delta_n$ and $\eta=f^{q_{n+1}}:I_n\to
\Delta_n$. This pair $(\xi, \eta)$ is what we call a {\it critical
commuting pair}. Each $f_{n+1}$ is by the definition the {\it
renormalization without rescaling} of $f_n$. Conjugating $f_n$ by the
affine map that takes the critical point $c$ to $0$ and $I_n$ to
$[0,1]$ we obtain ${\Cal R}^n(f)$, the $n$-th renormalization of $f$. 

\subhead
Holomorphic pairs
\endsubhead
The concept of {\it holomorphic commuting pair} was introduced in
{\refdF}, and will play a crucial role in this paper. We recall the
definition and some of the relevant facts about these objects,
henceforth called simply {\it holomorphic pairs}.
Assume we are given a configuration of four simply-connected domains
${\Cal O}_\xi,{\Cal O}_\eta,{\Cal O}_\nu, {\Cal V}$ in the complex
plane, called a {\it bowtie}, such that 

\itemitem{($a$)} Each ${\Cal O}_\gamma$ is a Jordan domain whose
closure is contained in ${\Cal V}$; 
\itemitem{($b$)} We have $\overline{{\Cal O}}_\xi\cap \overline{{\Cal
O}}_\eta= \{0\}\subseteq {\Cal O}_\nu$;
\itemitem{($c$)} The sets ${\Cal O}_\xi\setminus {\Cal O}_\nu$, ${\Cal
O}_\eta\setminus {\Cal O}_\nu$, ${\Cal O}_\nu\setminus {\Cal O}_\xi$
and ${\Cal O}_\nu\setminus {\Cal O}_\eta$ are non-empty and connected.

\noindent A holomorphic pair with domain ${\Cal U}={\Cal O}_\xi\cup
{\Cal O}_\eta\cup {\Cal O}_\nu$ is the dynamical system generated by
three holomorphic maps $\xi:{\Cal O}_\xi\to {\Bbb C}$, $\eta:{\Cal
O}_\eta\to {\Bbb C}$ and $\nu:{\Cal O}_\nu\to {\Bbb C}$ satisfying the
following conditions (see Figure 1).

\itemitem{[H$_1$]\enspace} Both $\xi$ and $\eta$ are univalent onto
${\Cal V}\cap {\Bbb C}(\xi(J_\xi))$ and ${\Cal V}\cap {\Bbb
C}(\eta(J_\eta))$ respectively, where $J_\xi={\Cal O}_\xi\cap {\Bbb R}$
and $J_\eta={\Cal O}_\eta\cap {\Bbb R}$. (Notation: ${\Bbb C}(I)=({\Bbb
C}\setminus {\Bbb R})\cup I$.)

\itemitem{[H$_2$]\enspace} The map $\nu$ is a $3$-fold branched cover
onto ${\Cal V}\cap {\Bbb C}(\nu(J_\nu))$, where $J_\nu={\Cal
O}_\nu\cap {\Bbb R}$, with a unique critical point at $0$.

\itemitem{[H$_3$]\enspace} We have ${\Cal O}_\xi\ni \eta(0)<0<\xi(0)\in
{\Cal O}_\eta$, and the restrictions $\xi|[\eta(0),0]$ and
$\eta|[0,\xi(0)]$ constitute a critical commuting pair.

\itemitem{[H$_4$]\enspace} Both $\xi$ and $\eta$ extend holomorphically
to a neighborhood of zero, and we have
$\xi\circ\eta(z)=\eta\circ\xi(z)=\nu(z)$ for all $z$ in that
neighborhood. 

\itemitem{[H$_5$]\enspace} There exists an integer $m\ge 1$, called the
{\it height} of $\Gamma$, such that $\xi^m(a)=\eta(0)$, where $a$ is
the left endpoint of $J_\xi$; moreover, $\eta(b)=\xi(0)$, where $b$ is
the right endpoint of $J_\eta$.

The interval $J=[a,b]$ is called the {\it long dynamical interval} of
$\Gamma$, whereas $\Delta=[\eta(0),\xi(0)]$ is the {\it short dynamical
interval} of $\Gamma$. They are both forward invariant under the
dynamics. The {\it rotation number} of $\Gamma$ is by definition the
rotation number of the critical commuting pair of $\Gamma$ (condition
H$_3$).

\bigskip
$$
\psannotate{
\psboxto(11cm;0cm){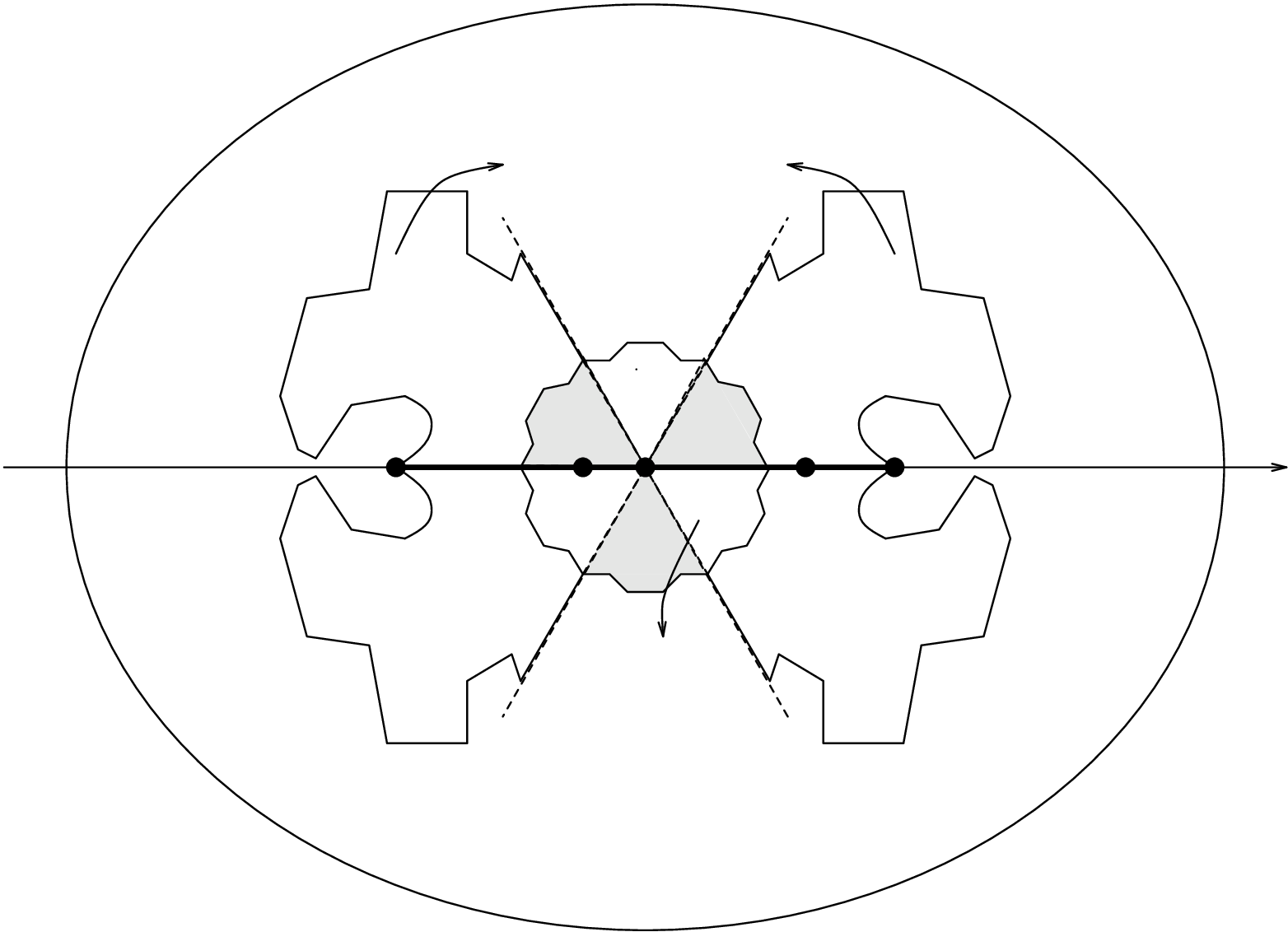}}
{\fillinggrid\at{7\pscm}{10\pscm}{${\Cal O}_\xi$}
             \at{20.6\pscm}{10\pscm}{${\Cal O}_\eta$}
             \at{10.5\pscm}{14\pscm}{$\xi$}
             \at{16.7\pscm}{14\pscm}{$\eta$}
             \at{20\pscm}{15\pscm}{${\Cal V}$}
             \at{12.6\pscm}{11\pscm}{${\Cal O}_\nu$}
             \at{12.7\pscm}{9\pscm}{$0$}
             \at{10.8\pscm}{7.5\pscm}{$\eta(\!0\!)$}
             \at{14.8\pscm}{7.5\pscm}{$\xi(\!0\!)$}
             \at{7.4\pscm}{7.5\pscm}{$a$}
             \at{16.2\pscm}{7.5\pscm}{$b$}
             \at{12.2\pscm}{5\pscm}{$\nu$}
}
$$
\smallskip
\centerline{\it Figure 1}
\bigskip

Examples of holomorphic pairs with arbitrary rotation number and
arbitrary height were carefully constructed in \S IV of {\refdF}, with
the help of the family of entire maps given by
$$
z\mapsto z+\theta-\frac{1}{2\pi}\sin{2\pi z}
$$
where $\theta$ is real, the so-called {\it Arnold family}. Many
interesting properties of holomorphic pairs can be proved with these
examples at hand. 

The {\it shadow} of a holomorphic pair is the map $F:{\Cal O}_\xi\cup
{\Cal O}_\eta\cup {\Cal O}_\nu\to {\Cal V}$ given by
$$
F(z)\;=\;
\cases
\xi(z),&\text{when $z\in {\Cal O}_\xi$}\cr
{}&{}\cr
\eta(z),&\text{when $z\in {\Cal O}_\eta$}\cr
{}&{}\cr 
\nu(z),&\text{when $z\in {\Cal O}_\nu\setminus ({\Cal O}_\xi\cup {\Cal
O}_\eta)$}\cr
\endcases
$$
The shadow captures the essential dynamical features of a holomorphic
pair. 
It is a simple fact (see {\refdF}, Prop.~II.4) that every
$\Gamma$-orbit is an $F$-orbit and conversely. 

The {\it limit set} of a holomorphic pair $\Gamma$, denoted ${\Cal
K}_\Gamma$, is the closure of the set of points in ${\Cal U}$ which
never escape from ${\Cal U}$ under iteration by $F$. The annulus ${\Cal
V}\setminus {\Cal U}$ is a fundamental domain for the dynamics off the limit
set, and for this reason it will be called the {\it fundamental annulus} of
$\Gamma$.

Holomorphic pairs can be renormalized. In other words, the first
renormalization of the critical commuting pair of $\Gamma$ extends in
a natural way to a holomorphic pair ${\Cal R}(\Gamma)$ with the same
co-domain ${\Cal V}$. For the careful construction of ${\Cal
R}(\Gamma)$, see Prop.~II.3 in {\refdF}.
There is also a pull-back theorem for holomorphic pairs. Let us say
that $\Gamma$ has {\it geometric boundaries} if $\partial {\Cal U}$
and $\partial {\Cal V}$ are quasicircles.

\proclaim{Theorem {2.3}} Let $\Gamma$ and $\Gamma'$ be holomorphic
pairs with geometric boundaries and let $h:J\to J'$ be a
quasisymmetric conjugacy between $F|J$ and $F'|J'$. Then there exists
a quasiconformal conjugacy $H:{\Cal V}\to {\Cal V}'$ between $\Gamma$ 
and $\Gamma'$ which is an extension of $h$.
\endproclaim

For a proof of this theorem, see Th.~III.1 in {\refdF}. 

\subhead
Conformal distortion
\endsubhead
We will need some classical facts about conformal maps, see {\refAh} and
{\refMca}. Consider a univalent map $\phi:\Omega\to {\Bbb C}$ on a domain
$\Omega\subseteq {\Bbb C}$. The {\it distortion} of $\phi$ on a
compact set $E\subseteq \Omega$ is the largest ratio
$|\phi'(z)|/|\phi'(w)|$ with $z,w\in E$. When $E$ is convex (a disk in
most cases), the distortion on $E$ is bounded by $\exp{N_\phi(E)}$,
where $N_\phi(E)=\sup_{z\in E}{|\phi''(z)|/|\phi'(z)|}$ is the {\it
total non-linearity} of $\phi$ on $E$. The {\it Koebe distortion
theorem} states that $N_\phi(E)\le 2/\dist (E,\partial\Omega )$. 
{\it Koebe's one-quarter theorem} states that if $D\subseteq \Omega$
is a disk centered at $z$, then $\phi(D)$ contains a disk centered at
$\phi(z)$ of radius $\frac{1}{4}\diam(D)|\phi'(z)|$.

\subhead
Almost parabolic maps
\endsubhead
In order to adapt Yampolsky's argument for the complex bounds in section 3, we
will need to know some general facts about complex analytic maps that are very
close to maps with a parabolic fixed-point.
Given $J\subseteq {\Bbb R}$ and $0<\theta<\pi$, we denote by
${\Bbb P}_\theta(J)\subseteq {\Bbb C}$ the set of all $z$ in
the complex plane that {\it view} $J$ under an angle $\ge\theta$. This set is 
the {\it Poincar\'e neighborhood} of $J$ with angle $\theta$.

\demo{Definition} Let $J\subseteq {\Bbb R}$ be an interval, and let
$\theta>0$. A holomorphic univalent map $\phi:{\Bbb P}_\theta(J)\to
{\Bbb C}$ is called {\it almost parabolic} if the following conditions
are satisfied. 
\itemitem{(a)} $\phi$ is symmetric about the real axis.
\itemitem{(b)} $\phi|J$ is monotone without fixed points.
\itemitem{(c)} $\phi$ has positive Schwarzian derivative on $J$.
\itemitem{(d)} $J\cap\phi(J)$ is non-empty.

\noindent If $\Delta_\phi$ is the interval $J\setminus\phi(J)$, the largest
$a=a(\phi)>0$ such that $\phi^{a-1}(\Delta_\phi)\subseteq J$ is called the
{\it length} of $\phi$. The number $\theta=\theta(\phi)$ is the {\it angle} of
$\phi$. 
\enddemo

Given $0<\sigma<1$, we denote by ${\Cal F}_\sigma$ the family of all
almost parabolic maps $\phi$ such that $|\Delta_\phi|\ge \sigma|J|$
and $|\phi^{a-1}(\Delta_\phi)|\ge \sigma|J|$, and also normalized so that
$$
[0,1]\;=\;\Delta_\phi\cup\phi(\Delta_\phi)\cup\cdots
\cup\phi^{a-1}(\Delta_\phi) 
\ .
$$
A fundamental lemma due to Yoccoz (see {\refdFdM} for a proof) states that for
each $0\le j \le a-1$ we have 
$$
\frac{1}{C_\sigma m(j)^2}\;\le\;|\phi^j(\Delta_\phi)|\;\le\;
\frac{C_\sigma}{m(j)^2} 
\ ,
$$ 
where $C_\sigma>1$ and $m(j)=\min\{j+1, a-j\}$. Every member of ${\Cal 
F}_\sigma$ whose length is sufficiently large has two fixed
points, symmetric about the real axis. More precisely, we have the following
fact. 

\proclaim{Lemma {2.4}} Given $0<\sigma<1$, there exist $C>1$, $a_0>0$ and
$\theta_0>0$ such that, if $\phi\in {\Cal F}_\sigma$ has length
$a=a(\phi)>a_0$ and angle $\theta(\phi)<\theta_0$, then 
there exist two attracting fixed points $z_+\in {\Bbb H}\cap
\dom(\phi)$ and $z_-=\overline{z}_+$ with 
$$
\frac{1}{Ca}\;\le\;\Im{z_+}\;\le\;\frac{C}{a}
\ .
$$
Moreover, if $|z-z_+|\le C/a$ then $|z-\phi(z)|\le C/a^2$.
\endproclaim

\demo{Proof} Follows from Yoccoz's lemma, the saddle-node bifurcation
and a normality argument. \qed
\enddemo

The family ${\Cal F}_\sigma$ is normal in the sense of Montel, and
every limit is a map with a parabolic (indifferent) fixed point on the
real axis.


\proclaim{Lemma {2.5}} Given $W\subseteq {\Bbb H}$ compact and an open set
$D\supseteq [0,1]$ in the plane, there exist $N_*>0$, $\theta_*>0$ and $a_*>0$
with the following property. For each $\phi\in 
{\Cal F}_\sigma$ such that $a(\phi)\ge a_*$ and $\theta(\phi)<\theta_*$, the
domain of $\phi$ contains $W$, and for each $z\in W$ there exists $n<N_*$
such that $\phi^n(z)\in D$. 
\endproclaim

\demo{Proof} If the statement is false, we find sequences $N_k\to\infty$,
$a_k\to\infty$ and $\theta_k\to 0$, maps $\phi_k\in{\Cal F}_\sigma$ with
$\theta(\phi_k)=\theta_k$ and $a(\phi_k)=a_k$ (whose domains contain $W$), and
points $z_k\in W$ such that $\phi_k^n(z_k)$, whenever defined, does not belong
to $D$ for all $n\le N_k$. Since ${\Cal F}_\sigma$ is normal and $W$ is
compact, we can assume that $\phi_k\to\phi: {\Bbb H}\to {\Bbb H}$ uniformly
on compacta and that $z_k\to z\in W$. Applying Lemma {2.4} to each $\phi_k$,
we deduce that $\phi$ has a fixed-point $x_0\in [0,1]$. By the Denjoy-Wolff
theorem, $\phi^n(z)\to x_0$ as $n\to\infty$. Hence there exists $N$ such that
$\phi^N(z)\in D$. But then $\phi_k^N(z_k)\in D$ also, for all sufficiently
large $k$, a contradiction.\qed
\enddemo

\heading 
3. Complex bounds for real-analytic circle maps
\endheading

We will extend Yampolsky's proof of complex bounds for circle
maps in the Epstein class {\refYam} to the more general case of {\it
real-analytic} circle maps. More precisely, we will
prove the following result.

\proclaim{Theorem {3.1}} (The complex bounds) Let $f: S^1\to S^1$ be
a real-analytic critical circle map with arbitrary irrational rotation
number. Then there exists $n_0=n_0(f)$ such that for all $n\ge n_0$ the 
$n$-th renormalization of $f$ extends to a holomorphic pair with geometric
boundaries whose fundamental annulus has conformal modulus bounded from below
by a universal constant.  
\endproclaim

The main step in the proof of Theorem {3.1} is to show that for all
sufficiently large $n$ the appropriate inverse branch of $f^{q_{n+1}}$
maps a sufficiently large disk around the $n$-th renormalization
domain $I_n\cup I_{n+1}$ well within itself.

We consider the unit circle $S^1={\Bbb R}/{\Bbb Z}$ embedded in the
infinite cylinder ${\Bbb C}/{\Bbb Z}$, and in all the geometric
considerations that follow we use on ${\Bbb C}/{\Bbb Z}$ the conformal
metric induced from the standard Euclidean metric $|dz|$ on ${\Bbb C}$
via the exponential map ${\Bbb C}\to {\Bbb C}/{\Bbb Z}$. Note that
$\Im {z}$ is well-defined for every $z\in {\Bbb C}/{\Bbb Z}$ (as the
imaginary part of any one of its pre-images under the exponential).

In what follows, we will fix $f:S^1\to S^1$ as in the statement of
Theorem {3.1}. Since $f$ is real-analytic, it extends to a
holomorphic map $f: A_R\to {\Bbb C}/{\Bbb Z}$, where $A_R$ is
the annulus $\{z\in {\Bbb C}/{\Bbb Z}: |\Im{z}| <R\}$. Making
$R$ smaller if necessary, we may assume that $f$ has no critical
points outside $S^1$. Using Koebe's distortion theorem, it is easy to
see that there exists $R_0>0$ such that, if $z\in S^1$ and $f(z)$ is
at a distance $> R_0$ from the critical value of $f$, then the
inverse branch $f^{-1}$ which maps $f(z)$ back to $z$ is well-defined
and univalent on the disk $D(f(z),R_0)$.

The key to the proof of Theorem {3.1} is to show that the $n$-th
renormalization of $f$ satisfies an inequality of the form
$|{\Cal R}^n(f)(z)|\ge C|z|^3$ on a neighborhood of the origin, where
$C$ is an absolute constant. The precise statement is the following.

\proclaim{Proposition {3.2}} There exist universal constants $B,C>0$
such that the following holds. For all $R>0$, there exists
$n_0=n_0(R)$ such that for all $n\ge n_0$ and all $z$ such that
$|z|>B$ and $|{\Cal R}^nf(z)|<R$ we have $|{\Cal R}^nf(z)|\ge C|z|^3$.
\endproclaim

This in turn depends on the following crucial proposition, which is
the analogue of Lemma 3.1 in {\refYam}.  
Following {\refYam}, for each $m\ge 1$ we define $D_m\subseteq {\Bbb
C}/{\Bbb Z}$ to be the disk of diameter the interval $[f^{q_{m+1}}(c), 
f^{q_m-q_{m+1}}(c)]\subseteq S^1$ containing the critical point $c$.
Note that $\roman{diam}\,(D_m)$ is comparable with $|I_m|$.

\proclaim{Proposition {3.3}} There exist universal constants $B_1$ and
$B_2$ and for each $N\ge 1$ there exists $n(N)$  such that
for all $n\ge n(N)$ the inverse branch $f^{-q_{n+1}+1}$ taking
$f^{q_{n+1}}(I_n)$ back to $f(I_n)$ is well-defined and univalent over
$\Omega_{n,N}=(D_{n-N}\setminus S^1)\cup f^{q_{n+1}}(I_n)$, and for
all $z\in \Omega_{n,N}$ we have   
$$
\frac{\roman{dist}\,\left(f^{-q_{n+1}+1}(z), f(I_n)\right)}{|f(I_n)|}
\;\le\; B_1\left(\frac{\roman{dist}\,(z,I_n)}{|I_n|}\right)+B_2
\ ,
\slabel{lgrowth}
$$
\endproclaim

Yampolsky's proof of this fact in {\refYam} is greatly
facilitated by his assumption that the circle map has the Epstein
property (namely, its lift to ${\Bbb C}$ has univalent inverse
branches globally defined in the upper half-plane). Thus, the
Poincar\'e neighborhoods ${\Bbb P}_\theta (J)$ are automatically
invariant in the sense that $f^{-1}({\Bbb P}_\theta(J))\subseteq {\Bbb
P}_\theta(f^{-1}J)$. In our setting this is no longer true, but as we
look at smaller and smaller scales we see that it is asymptotically
true, as the following lemma shows. 

\proclaim{Lemma {2.4}} For every small $a>0$, there exists
$\theta(a)>0$ satisfying $\theta(a)\to 0$ and $a/\theta(a)\to 0$ as
$a\to 0$, such that the following holds. Let $F:{\Bbb D}\to {\Bbb C}$
be univalent and symmetric about the real axis, and assume $F(0)=0$,
$F(a)=a$. Then for all $\theta\ge \theta(a)$ we have 
$F\left({\Bbb P}_\theta ([0,a])\right)\subseteq
{\Bbb P}_{(1-a^{1+\delta})\theta}([0,a])$, where $0<\delta<1$ is an
absolute constant. 
\endproclaim
\demo{Proof} 
There exists a Moebius transformation $G$ such that $G(0)=0$, $G(a)=a$
and $|D^jF(x)-D^jG(x)|\le C_0a^{3-j}$ for $j=0,1,2$ and all $x\in
[0,a]$, where $C_0$ is an absolute constant (consider the Moebius
transformation with the same $2$-jet as $F$ at zero, post-composed
with a linear map to meet the normalization condition $G(a)=a$). This
$G$ has no pole in a disk $D_0\subseteq {\Bbb D}$ of definite
radius around zero.
Let $\varphi(z)=F(z)-G(z)$. Then $\varphi(z)=b_1z+b_2z^2+\cdots$ for
all $z\in D_0$, where $|b_1|\le C_0a^2$ and $|b_2|\le C_0a$. Take a
small number $\varepsilon>0$ and consider the disk
$D_1=D(0,a^{1-\varepsilon})\subseteq D_0$. In this disk we have the
estimate 
$$
|\varphi(z)|\;\le\;|b_1||z|+|b_2||z|^2+C_1|z|^3\;\le\;
C_2a^{2-2\varepsilon}|z|
$$
At the same time, using the fact that $F'(\zeta)=1$ for some $\zeta\in
[0,a]$ and the Koebe distortion lemma, we see that $|F(z)|\ge C_3|z|$
and $|G(z)|\ge C_3|z|$ for all $z\in D_1$. Therefore, for every such
$z$ the triangle with vertices at $0$, $F(z)$ and $G(z)$ has an angle
at zero $\le C_4a^{2-2\varepsilon}$. Similarly for the angle at $a$ in
the triangle with vertices $a$, $F(z)$ and $G(z)$. Now suppose $z\in
{\Bbb P}_\theta([0,a])$, for $\theta\ge a^\varepsilon$. Since $G$
preserves this neighborhood, $G(z)$ forms an angle $\ge \theta$ with
$(-\infty,0]$, and the same holds for the angle $G(z)$ forms with
$[a,+\infty)$. It follows that $F(z)\in {\Bbb P}_{\theta'}([0,a])$,
where $\theta'=\theta-C_4a^{2-2\varepsilon}\ge \theta(1-a^{1+\delta})$
for some $0<\delta<1$ depending only on $\varepsilon$. This proves the
lemma with $\theta(a)=a^{\varepsilon}$.\qed
\enddemo

The main consequence of this lemma is the following.

\proclaim{Lemma {2.5}} For each $n\ge 1$ there exist $K_n\ge 1$ and
$\theta_n>0$, with $K_n\to 1$ and $\theta_n\to 0$ as $n\to\infty$,
such that for all $\theta\ge\theta_n$ and all $1\le j\le q_{n+1}$
the inverse branch $f^{-j+1}$ mapping $f^{j}(I_n)$ back to $f(I_n)$ is
well-defined over ${\Bbb P}_\theta(f^{j}(I_n))$ and maps this
neighborhood univalently into ${\Bbb P}_{\theta/K_n}(f(I_n))$.
\endproclaim

\demo{Proof} Let $d_n=\max_{1\le j\le q_{n+1}}{|f^j(I_n)|}$; from the
real bounds, these numbers go to zero exponentially with $n$. Take
$\delta >0$ as in Lemma {2.4}, and let $K_n$ be given by
$$
K_n^{-1}\;=\;\prod_{j=1}^{q_{n+1}}\left(1-|f^j(I_n)|^{1+\delta}\right) 
\ .
$$
Then define $\theta_n=K_n\theta(d_n)$, where $\theta(\cdot)$ is the
function in Lemma {2.4}.
Note that
$$
\log{K_n}\;\le\;C\,\sum_{j=1}^{q_{n+1}}|f^j(I_n)|^{1+\delta}
\;\le\;Cd_n^{\delta} 
\ .
$$
Therefore $K_n\to 1$ and $\theta_n\to 0$ as required. Also,
$d_n/\theta_n\to 0$.

Now fix $j$ as in the statement and suppose $\theta\ge \theta_n$.
Define inductively $\vartheta_0=\theta$ and
$\vartheta_{i+1}=(1-|f^{j-i}(I_n)|^{1+\delta})\vartheta_i$ for $i=0,
1,\ldots , j-2$, and note that $\vartheta_{j-1}\ge \theta/K_n$.
Moreover,
$$
\roman{diam}\,\left({\Bbb
P}_{\vartheta_i}\left(f^{j-i}(I_n)\right)\right)\;=\;  
\frac{|f^{j-i}(I_n)|}{\sin{\vartheta_i}}\;\le\;
\frac{C'd_n}{\theta(d_n)}\;<\!<\;R_0
\ .
$$
Therefore $f^{-1}$ is well-defined and univalent over ${\Bbb
P}_{\vartheta_i} (f^{j-i}(I_n))$, and by Lemma {2.4} we have the
inclusion  
$f^{-1}({\Bbb P}_{\vartheta_i}(f^{j-i}(I_n)))\subseteq
{\Bbb P}_{\vartheta_{i+1}}(f^{j-i-1}(I_n))$. This completes the
proof.\qed 

\enddemo

\demo{Remark} The same result holds true if we replace $I_n$ by any
interval $J\supseteq I_n$ such that the map $f^{q_{n+1}-1}:f(J)\to
f^{q_{n+1}}(J)$ is a diffeomorphism. 
\enddemo

We will need four lemmas concerning the sequence $\{D_m\}$ introduced
earlier. 
The first is an easy consequence of Lemma {2.5} and the above
remark. 

\proclaim{Lemma {2.6}} There exists $m_0\ge 1$ such that for all
$m\ge m_0$ the inverse branch $f^{-q_m+1}$ taking $f^{q_m}(I_m)$ back
to $f(I_m)$ is well-defined and univalent in $D_m$, and 
$$
\frac{\diam(f^{-q_m+1}(D_m))}{|I_m|}\;\le\;
C\,\frac{\diam(f(D_m))}{|f(I_m)|}
\ .
$$
\endproclaim

The second is the analogue of Lemma 4.1 in {\refYam}.

\proclaim{Lemma {2.7}} There exist $\varepsilon_1>0$ and $m_1\ge
m_0$ such that for all $m\ge m_1$ and each $w\in
f^{-q_{m+1}}(D_m)\setminus D_m$ we have ($a$)
${\roman{dist}}\,(w,I_m)\le C|I_m|$ and ($b$) for each $x\in I_m$,
$\varepsilon_1 <|\arg{(w-x)}|<\pi-\varepsilon_1$. 
\endproclaim
\demo{Proof} The same proof given in {\refYam} applies here.
Invariance of Poincar\'e neighborhoods is replaced by {\it
quasi-invariance}, using Lemma {2.5}. \qed
\enddemo
\bigskip
\bigskip
$$
\psannotate{
\psboxto(10cm;0cm){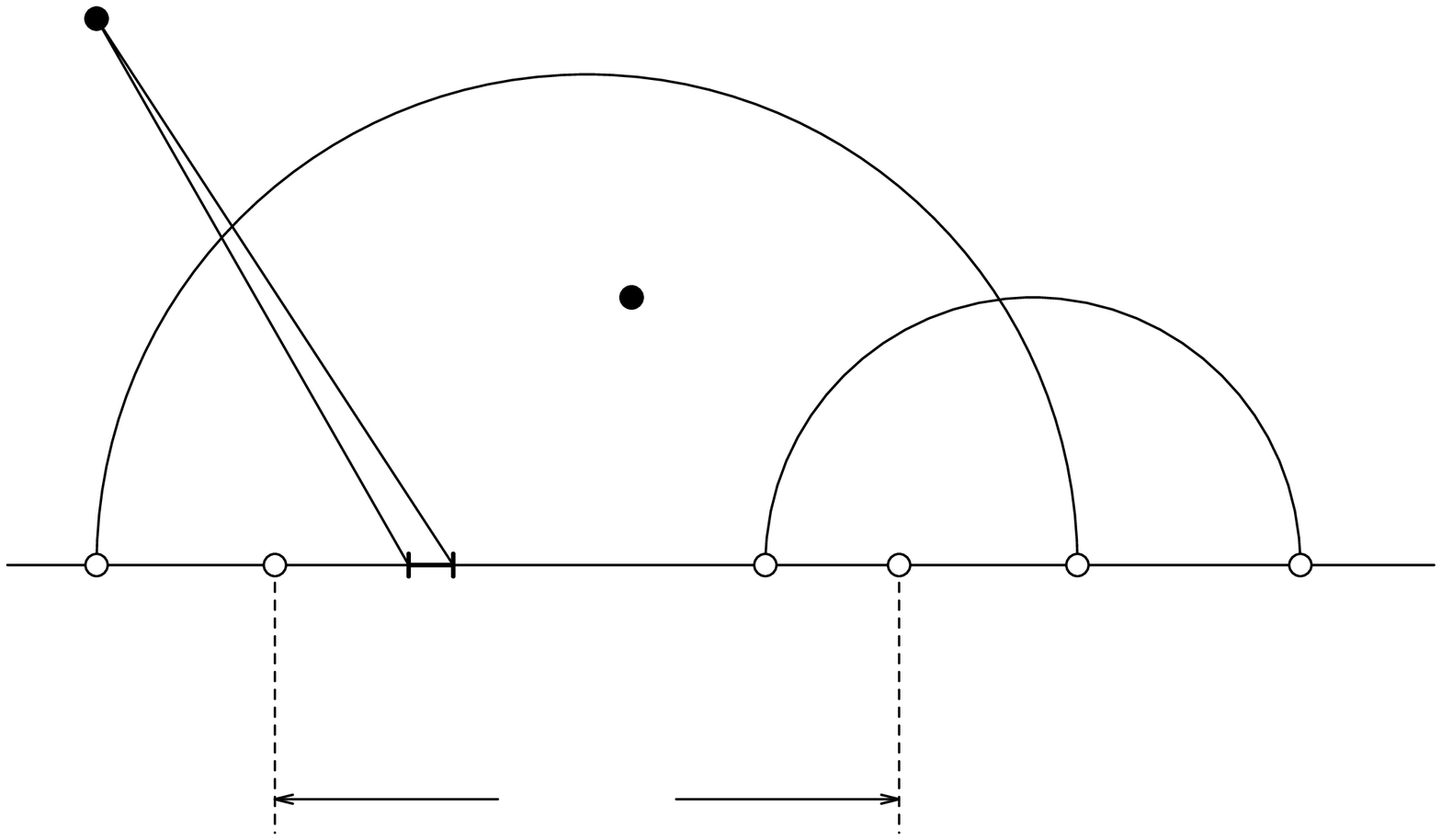}}
{\fillinggrid\at{10.15\pscm}{0.3\pscm}{$I_m$}
             \at{1.1\pscm}{2.8\pscm}{$f^{q_m-q_{m+1}}(c)$}
             \at{12.7\pscm}{2.8\pscm}{$I_{m+2}$}
             \at{14.95\pscm}{2.8\pscm}{$I_{m+1}$}
             \at{17.9\pscm}{2.8\pscm}{$f^{q_{m+1}-q_{m+2}}(c)$}
             \at{13.7\pscm}{4.3\pscm}{$c$}
             \at{2.75\pscm}{11.45\pscm}{$\zeta''$}
             \at{9.3\pscm}{11.1\pscm}{$D_m$}
             \at{9.8\pscm}{7.5\pscm}{$\zeta$}
             \at{17.1\pscm}{7.1\pscm}{$D_{m+1}$}
}
$$
\smallskip
\centerline{\it Figure 2}
\bigskip

The third is the analogue of Lemma 4.4 in {\refYam}. It provides
us with the tools we need for the inductive step in the proof of
Proposition {3.3} (see Figure 2). 

\proclaim{Lemma {2.8}} There exist $\varepsilon_2>0$ and $m_2\ge
m_0$ such that the following holds for all $m\ge m_2$. Let $\zeta\in
D_m\setminus D_{m+1}$ be a point not on the circle, and let
$\zeta'=f^{-q_m}(\zeta)$ and $\zeta''=f^{-q_{m+2}}(\zeta')$. Then
either $\zeta''\in D_{m+1}$, or else
${\roman{dist}}\,(\zeta'',I_{m+1})\le C|I_m|$ and
$\varepsilon_2<\arg{(\zeta''-x)}<\pi-\varepsilon_2$ for all $x\in
I_m\cup I_{m+1}$.
\endproclaim
\demo{Proof} Once again, the proof of Lemma 4.4 in {\refYam} can be
repeated here, {\it mutatis mutandis}.
\qed
\enddemo

\demo{Notation} Given a point $\zeta \in {\Bbb C}$ and an interval
$J=(a,b)\subseteq {\Bbb R}$, we denote by $\roman{angle}\,(\zeta,J)$
the smallest of the angles $\pi-\arg{(\zeta-a)}$ and
$\arg{(\zeta-b)}$.  
\enddemo

The fourth is a consequence of Lemma {2.5}.

\proclaim{Lemma {2.9}} There exist universal constants $N_*>0$ and $a_*>0$
and some $m_3>0$ with the following property. For all $m\ge m_3$ such that
$a_{m+1}>a_*$ and each $w\in V_m=f^{-q_{m+1}}(D_m)\setminus D_m$, there exists
$1<i<N_*$ such that the iterate $(f^{-q_{m+1}})^i(w)$ is well-defined and
belongs to $D_m$. 
\endproclaim

\demo{Proof} Normalize $f^{-q_{m+1}}$ so that $I_m\setminus I_{m+2}$ becomes
the interval $[0,1]$ to get an almost parabolic map $\phi_m$. Note that
$\phi_m\in {\Cal F}_\sigma$ for some $\sigma$ depending only on the real
bounds. Let $W_m\subseteq {\Bbb H}$ be the image of $V_m$ under such
normalization. By Lemma {2.7}, there exists a fixed compact set $W\subseteq
{\Bbb H}$ such that $W_m\subseteq W$ for all sufficiently large
$m$. Similarly, the normalized copies of $D_m$ contain a fixed open set
$D\supseteq [0,1]$ for all sufficiently large $m$. Hence we can take $N_*$ and
$a_*$ as given by Lemma {2.5}.\qed
\enddemo

\demo{Proof of Proposition {3.3}} 
We will start with a point $z$ in the disk $D_{n-N}$. For the argument 
to work, $n$ will have to be sufficiently large. We start taking
$n>N+\max\{m_1,m_2,m_3\}$, where $m_1$, $m_2$ and $m_3$ are given respectively
by Lemma {2.7}, Lemma {2.8} and Lemma {2.9}.

In keeping with Yampolsky's notation in {\refYam}, let us denote by
$J_{-i}$ the interval $f^{q_{n+1}-i}(I_n)$. 
Also, given $z$, let $z_{-i}=f^{-i}(z)$ be the corresponding
pre-images of $z$. 

The proof runs by finite induction in the range $n-N\le m\le
n$. Let $m$ be the largest with the property that $z\in D_m$, and keep
in mind that $\roman{dist}\,(z,I_n)\asymp |I_m|$. Consider
those moments $i_1<i_2<\cdots <i_\ell$ in the backward orbit
$\{J_{-i}\}$ before the first return to $I_{m+1}$ such that
$J_{-i_k}\subseteq I_m$.  Then, there are two possibilities. 

The first possibility is that $z_{-i_\ell}\notin D_m$. In this case there
exists a smallest $k\le \ell$ such that $z_{-i_s}\notin D_m$ for
$s=k,k+1,\ldots,\ell$. 
We claim that $|J_{-i_k}|\asymp |I_n|$. This is clear from the real bounds if
$\ell=a_{m+1}\le a_*$, where $a_*$ is given by Lemma {2.9}. If on the
other hand $\ell>a_*$, then again by Lemma {2.9} we must have
$\ell-k<N_*$, and the claim follows from Lemma {2.2}. Therefore, by Lemma
{2.7}, 
$$
\frac{\roman{dist}\,(z_{-i_k},J_{-i_k})}{|J_{-i_k}|}\;\le\; 
C\frac{|I_m|}{|J_{-i_k}|}\;\le\;C'\frac{\roman{dist}\,(z,I_n)}{|I_n|} 
\ .
\slabel{IK}
$$
Moreover, $\roman{angle}\,(z_{-i_k}, J_{-i_k})\ge \varepsilon_1$, so
there exists $\theta=\theta(\varepsilon_1,N)$ such that $z_{-i_k}\in
{\Bbb P}_\theta(J_{-i_k})$. Now, if $n$ is sufficiently large,
$\theta_n<\theta$ and we can use Lemma {2.5} to get that
$z_{-q_{n+1}+1}\in {\Bbb P}_{\theta/K_n}(f(I_n))$. This gives us
$$
\frac{\roman{dist}\,(z_{-q_{n+1}+1},f(I_n))}{|f(I_n)|}\;\le\; 
C''K_n\frac{\roman{dist}\,(z_{-i_k},J_{-i_k})}{|J_{-i_k}|}
\ ,
\slabel{IKo}
$$
and this together with {\eqnoIK} yields the Proposition in this case.

The second possibility is that $\zeta=z_{-i_\ell}\in D_m$, and we can
assume that $\zeta\notin D_{m+1}$ (otherwise the induction step is
complete). In this case, consider $\zeta'=f^{-q_m}(\zeta)$ and
$\zeta''=f^{-q_{m+2}}(\zeta')$ and the corresponding interval
$J''=f^{-q_m-q_{m+2}} (J_{-i_\ell})$, and apply Lemma {2.8}. Then
either $\zeta''\in D_{m+1}$, in which case the induction step is complete, or
else $\roman{dist}\,(\zeta'',I_{m+1})\le C|I_m|$ and
$\roman{angle}\,(\zeta'', J'')\ge \varepsilon_1$, in which case we can
apply the same argument leading to {\eqnoIK} and {\eqnoIKo}. 

If the backward orbit survives all the steps of the induction, this
means that in the end $z_{-q_{n+1}+q_{n-1}}\in D_{n-1}$. By Lemma
{2.6}, the image of $D_{n-1}$  under $f^{-q_{n-1}+1}$ has diameter
comparable to $|f(I_n)|$, so the first member of {\eqnoIKo} is simply
bounded by an absolute constant. So in any case we have
{\eqnolgrowth}. \qed 
\enddemo

\heading 
4. Geometric control of holomorphic pairs
\endheading

In this section, we show that if a holomorphic pair satisfies certain
geometric constraints, then one can bound from below how much
its dynamics expands the hyperbolic metric in the complement of the
real axis. From the renormalization viewpoint, these constraints are
quite natural, being a consequence of the complex bounds. 

We start with the definition of geometric control. For brevity, let us
say that a finite set in the plane is $K$-bounded ($K>1$), if any two
non-zero distances $d_1$ and $d_2$ between points in the set satisfy
$d_1\le Kd_2$. We denote by $Q$ the map $z\mapsto z^3$ in the complex
plane. 

\demo{Definition} A holomorphic pair $\Gamma$ is controlled by
$K>1$, or simply $K$-controlled, if the following geometric conditions
are satisfied.
\itemitem{[G$_1$]} The set
$\left\{0,\xi(0),\eta(0),\nu(0),a,b\right\}$ 
is $K$-bounded.
\itemitem{[G$_2$]} We have $\diam ({\Cal V})\le K|J|$.
\itemitem{[G$_3$]} On the real axis, $|\gamma'(x)|\le K$ for $\gamma
=\xi,\eta,\nu$. 
\itemitem{[G$_4$]} $\xi$ has holomorphic extensions
$\xi_a:D(a,K^{-1}|J|)\to {\Bbb C}$ and $\xi_0:D(0,K^{-1}|J|)\to
{\Bbb C}$ such that $\xi_a=\psi_a\circ Q\circ\phi_a$ and
$\xi_0=\psi_\xi\circ Q$ respectively, where $\psi_a$, $\phi_a$ and
$\psi_\xi$ are univalent maps with distortion bounded by $K$, and
$\phi_a(a)=0$.
\itemitem{[G$_5$]} $\eta$ has holomorphic extensions
$\eta_b:D(b,K^{-1}|J|)\to {\Bbb C}$ and $\eta_0:D(0,K^{-1}|J|)\to
{\Bbb C}$ such that $\eta_b=\psi_b\circ Q\circ\phi_b$ and
$\eta_0=\psi_\eta\circ Q$ respectively, where $\psi_b$, $\phi_b$ and
$\psi_\eta$ are univalent maps with distortion bounded by $K$, and
$\phi_b(b)=0$.
\itemitem{[G$_6$]} The open sets ${\Cal O}_\xi\cap D(a,K^{-1}|J|)$ and
${\Cal O}_\eta\cap D(b,K^{-1}|J|)$ are connected.
\itemitem{[G$_7$]} We have $D(0,K^{-1}|J|)\subseteq {\Cal O}_\nu$.
\itemitem{[G$_8$]} We have $\mod ({\Cal V}\setminus {\Cal U})\ge K^{-1}$.
\enddemo
\demo{Remark} Condition G$_6$ guarantees that the shadow of a
$K$-controlled holomorphic pair extends to a single-valued map 
$$
F_{*}:\dom (\Gamma) \cup D(a,K^{-1}|J|)\cup D(b,K^{-1}|J|)\to {\Bbb C}
\ .
$$ 
We call this map the {\it extended shadow\/} of our holomorphic pair. 
\enddemo

Now a strengthened version of the complex bounds can be stated as
follows. 

\proclaim{Theorem {4.1}} There exists a universal constant $K_0>1$ such that
the following holds. If $\Gamma$ is a holomorphic
pair with arbitrary rotation number, there
exists $n_0=n_0(\Gamma)$ such that for all $n\ge n_0$ the holomorphic
pair ${\Cal R}^n\Gamma$ restricts to a holomorphic pair $\Gamma_n$ 
controlled by $K_0$. \qed
\endproclaim

Let us now formulate the tools from function theory that we need to
prove expansion of the hyperbolic metric by a controlled holomorphic
pair. 
Given a hyperbolic Riemann surface $X$, we write $\rho_X$ for the
hyperbolic density of $X$,  $d_X$ for the 
hyperbolic metric of $X$, and $\diam_X(E)$ for the diameter of
$E\subseteq X$ in that metric. If $f:{\Cal O}\subseteq X\to Y$ is a
holomorphic map between two hyperbolic Riemann surfaces (${\Cal O}$ an
open subset of $X$), we denote by 
$$
\|f'(z)\|_{X,Y}=\frac{\rho_Y(f(z))}{\rho_X(z)}\,|f'(z)|
$$
the derivative of $f$ at $z\in X$ measured with
respect to the hyperbolic metrics in $X$ and $Y$, and simplify the
notation a bit to $\|f'(z)\|_Y$ when $X=Y$. 

The key to expansion is the following lemma, stated as Prop.~4.9 in
{\refMcb}. 

\proclaim{Lemma {4.2}} There exists a positive function $C(s)$
decreasing to zero as $s$ decreases to zero such that the following
holds. If $f:X\subseteq Y$ is the inclusion of a hyperbolic 
Riemann surface into another, then $\|f'(z)\|_{X,Y}\le
C(d_Y(z,Y\setminus X))<1$.\qed
\endproclaim

The second tool we need, to be used in combination with the above, is
the following.

\proclaim{Lemma {4.3}} Let $\Omega$ be a doubly-connected region in the
plane, let $E$ be the bounded component of ${\Bbb C}\setminus\Omega$, and let
$V=E\cup \Omega$. Then we have
$$
\diam_V(E)\le\frac{C}{\mod (\Omega )}\ ,
\slabel{one}
$$
for some constant $C>0$, as well as 
$$
\mod (\Omega )\ge \frac{4}{\pi}\,\left[\frac{\delta}{\diam
(\Omega)}\right]^2\ ,
\slabel{two}
$$
where $\delta$ is the Euclidean distance between the boundary
components of $\Omega$.
\endproclaim
\demo{Proof} For {\eqnoone}, see {\refMca}, Th.~2.4, where in fact a
more precise estimate is given . For {\eqnotwo}, use the definition of
$\mod (\Omega)$ as the extremal length of the family ${\Cal F}$ of
curves in $\Omega$ joining the two components of $\partial\Omega$
to get 
$$
\mod(\Omega)\ge \frac{1}{A(\Omega)}\,\inf_{\gamma\in{\Cal
    F}}L^2(\gamma)
\ ,
$$
where $L(\gamma)$ is the Euclidean length of $\gamma$ and $A(\Omega)$
is the Euclidean area of $\Omega$. Since $A(\Omega)\le \pi (\diam
\Omega )^2/4$, whereas $L(\gamma)\ge \delta$ for all $\gamma\in{\Cal
  F}$, this proves {\eqnotwo}.\qed
\enddemo

While inequality {\eqnotwo} gives a lower bound for the modulus of an
annulus in terms of the minimum separation between boundary
components, an upper bound is provided by Teichm\"uller's inequality 
$$
\mod (\Omega)\le \Phi\left(\frac{\delta}{\diam\Omega}\right)
\ ,
\slabel{teich}
$$
where $\Phi$ is a certain universal monotone increasing function, {\it
  cf.} {\refAh}.

Later in this section it will also be very important to be able to
control the variation of holomorphic distortion of an analytic map
between hyperbolic Riemann surfaces. When the time comes, we will
resort to the following consequence of Koebe's distortion theorem, as
stated by McMullen in {\refMca}, Cor.~2.27.

\proclaim{Lemma {4.4}} Let $f:X\to Y$ be an analytic map between
hyperbolic Riemann surfaces whose derivative vanishes nowhere. Then
for all $x_1,x_2\in X$ we have
$$
\|f'(x_1)\|_{X,Y}^{\alpha}\le \|f'(x_2)\|_{X,Y}\le
\|f'(x_1)\|_{X,Y}^{1/\alpha} 
\ ,
$$
where $\alpha=\exp\{C_0d_X(x_1,x_2)\}$ and $C_0>0$ is a universal
constant. \qed
\endproclaim

Besides these very general tools, we will need the following specific
consequence of geometric control. Given $E\subseteq {\Bbb C}$ and
$\varepsilon>0$, we write $V(E,\varepsilon )
=\left\{z:d(z,E)<\varepsilon\,\diam(E)\right\}$. 

\proclaim{Lemma {4.5}} Given $K>1$, there exist $\varepsilon_0>0$
and $M>1$ such that the following holds for all
$0<\varepsilon\le\varepsilon_0$. If $\Gamma$ is a $K$-controlled
holomorphic pair and $F_*$ is its extended shadow, then 
\itemitem{$(a)$} $V(J,\varepsilon )\subseteq \dom (F_*)\/$;
\itemitem{$(b)$} For all $z\in V(J,\varepsilon )$ we have
$|F_*'(z)|\le M\/$;
\itemitem{$(c)$} There exists $0<\tau=\tau(\varepsilon,
K)<\varepsilon$ such that $F_*(V(J,\tau ))\subseteq
V(F_*(J),\varepsilon)\/$. 
\itemitem{$(d)$} $F_*(V(J,\varepsilon ))\subseteq \dom (\Gamma )\/$.
\endproclaim

\demo{Proof} Let $J_*=F_*(J)=[\xi(a),\eta(b)]=[\xi(a),\xi(0)]$.

To prove ($a$), we note that 
$$
\mod ({\Cal V}\setminus J_*)>\mod ({\Cal V}\setminus\dom (\Gamma ))
\ge K^{-1} 
\ ,
$$
by condition G$_8$. Hence, there exists
$\varepsilon_1>0$ depending only on $K$ such that
$V(J_*,\varepsilon_1)\subseteq {\Cal V}$, by Teichm\"uller's
inequality {\eqnoteich} and condition G$_2$. Using condition G$_4$ and
making $\varepsilon_1$ smaller if necessary, we find
$0<\varepsilon_2<K^{-1}$ and $0<\varepsilon_3<\varepsilon_1$, both
depending only on $K$, such that 
$$
\cases
D(\xi(a),\varepsilon_3|J_*|)\subseteq\xi_a(D(a,\varepsilon_2|J|))\subseteq
V(J_*,\varepsilon_1)&{}\cr
{}&{}\cr
D(\xi(0),\varepsilon_3|J_*|)\subseteq\xi_0(D(0,\varepsilon_2|J|))\subseteq
V(J_*,\varepsilon_1)&{}\cr
\endcases
\ .
\slabel{disks}
$$

Now, consider $a<x<0$. If $|\xi(x)-\xi(a)|<\varepsilon_3|J_*|$, then
by {\eqnodisks} we have $x\in D(a,\varepsilon_2|J|)$. Similarly, if
$|\xi(x)-\xi(0)|<\varepsilon_3|J_*|$, then $x\in
D(0,\varepsilon_2|J|)$. If neither of these happen, then 
$$
\xi(a)+\varepsilon_3|J_*|<\xi(x)<\xi(0)-\varepsilon_3|J_*|
\ .
$$
In this case the disk $D_x=D(\xi(x),\varepsilon_3|J_*|)\subseteq
V(J_*,\varepsilon_1)$  is contained in the image of $\xi$. Since
$|\xi'(x)|\le K$ by condition G$_3$, and since $|J_*|\ge K^{-1}|J|$ by
condition G$_1$, we deduce from Koebe's one-quarter theorem that
$$
\xi^{-1}(D_x)\supseteq D\left(x,\frac{\varepsilon_3|J|}{4K^2}\right)
\ .
$$
Putting these facts together, we get
$$
\xi^{-1}(V(J_*,\varepsilon_1))\supseteq V([a,0],\varepsilon')
\ ,
\slabel{xiv}
$$
where
$\varepsilon'=\min\{\varepsilon_2,\varepsilon_3/4K^2\}$. Proceeding 
similarly with $\eta$ replacing $\xi$ and condition G$_5$ replacing
condition G$_4$, and noting that $J_*\supseteq [\eta(0),\eta(b)]$, we
obtain $\varepsilon''>0$ depending only on $K$ such that 
$$
\eta^{-1}(V(J_*,\varepsilon_1))\supseteq V([0,b],\varepsilon'')
\ .
\slabel{etav}
$$
Taking $\varepsilon_0=\min\{|a|\varepsilon',|b|\varepsilon''\}/|J|$,
we deduce from {\eqnoxiv} and {\eqnoetav} that
$V(J,\varepsilon_0)\subseteq \dom (F_*)$, and this proves ($a$).

The proof of the remaining assertions is more of the same, so we omit
the details. Using condition G$_3$ and Koebe distortion, we see that
the maps defining $F_*$ have distortion on $V(J,\varepsilon_0/2)$
bounded only in terms of $K$. Therefore ($b$) holds for some constant 
$M$ depending only on $K$ if we replace $\varepsilon_0$ by
$\varepsilon_0/2$. Making this new $\varepsilon_0$ still smaller if
necessary, we can assume also that $V(J_*,\varepsilon_0)\subseteq\dom
(\Gamma)$.
Then, ($c$) follows for some $\tau$ as stated if we use ($b$) and
apply the mean value theorem to $F_*$ in $V(J,\tau)$. Finally,
replacing $\varepsilon_0$ by $\tau(\varepsilon_0,K)$, we see that
($d$) follows from ($c$).\qed
\enddemo

Next, let $\Gamma$ be a holomorphic pair and let $Y={\Cal V}\setminus
{\Bbb R}$. How do we compare the hyperbolic metric of $Y$ with the
hyperbolic metric of the upper half-plane ${\Bbb H}$?
The answer is given by the following lemma. 

\proclaim{Lemma {4.6}} Let $\rho_Y(z)$ and
$\rho(z)=1/|\Im{z}|$ be the hyperbolic densities of $Y$ and ${\Bbb
C}\setminus {\Bbb R}$, respectively. Then for all $z\in \dom (\Gamma)$ we
have  $\rho(z)<\rho_Y(z)<c_M\rho(z)$, where $0<c_M< 1$ is a constant
depending only on $M=\mod ({\Cal V}\setminus\dom (\Gamma))$.
\endproclaim

\demo{Proof} 
By monotonicity of hyperbolic densities, $\rho_Y>\rho$.
Let $\delta$ be the Euclidean
distance between $\partial {\Cal V}$ and $\partial\dom (\Gamma)$. By
Teichm\"uller's inequality {\eqnoteich}, we have 
$$
\delta\ge\Phi^{-1}(M)\,\diam {\Cal V}\ge \Phi^{-1}(M)\,|\Im{z}|
\ .
$$
Hence the largest disk $D$ contained in $Y$ and centered at $z$
has radius 
$$
R\ge\min \{\delta , |\Im{z}|\}\ge 2c_M^{-1}|\Im{z}|
\ ,
$$
where $c_M^{-1}=\frac{1}{2}\,\min\{1,\Phi^{-1}(M)\}$. Therefore, again
by monotonicity of hyperbolic densities, we get
$$
\rho_Y(z)\le \rho_D(z)=\frac{2}{R}\le c_M\,\frac{1}{|\Im{z}|}=
c_M\,\rho(z)
\ .
$$
\enddemo

We are ready for the first expansion result.

\proclaim{Proposition {4.7}} Let $\Gamma$ be a $K$-controlled holomorphic
pair with co-domain ${\Cal V}$ and let $Y={\Cal V}\setminus{\Bbb R}$. Then for
each $\varepsilon>0$ we have the following.
\itemitem{$(a)$} There exists $\varrho=\varrho(\varepsilon,K)>0$
satisfying $\|F'(z)\|_Y\ge 1+\varrho$ for all $z$ in the domain of $\Gamma$
such that $|\Im{z}|\ge\varepsilon |J|$. 
\itemitem{$(b)$} There exists $\lambda=\lambda (\varepsilon,K)>1$ such that
$\|F'(z)\|_Y\ge \lambda$ for all $z\in \dom (\Gamma )\cap V(J,\varepsilon)$
such that $F(z)\notin V(J,\varepsilon)$. 

\noindent
Moreover, if $\Gamma$ is the renormalization of a holomorphic pair $\Gamma^*$
with co-domain ${\Cal V}^*$ and $Y^*={\Cal V}\setminus {\Bbb R}$, then the
same statements are true with $\|\cdot\|_{Y^*}$ replacing $\|\cdot\|_Y$
throughout.  
\endproclaim

\demo{Proof} It suffices to consider only
$\varepsilon\le\varepsilon_0$, where $\varepsilon_0$ is given by Lemma
{4.5}. By condition G$_8$ and (3), there exists
$\varepsilon_1=\varepsilon_1(K)>0$ such that the Euclidean distance
between $\partial\dom (\Gamma )$ and $\partial{\Cal V}$ is at least
$\varepsilon_1|J|$.
 
Let $E$ be the set of points in $\dom(\Gamma )^+\subseteq {\Cal V}^+$
such that $\Im{z}\ge\varepsilon |J|$, and let $\delta$ be the
Euclidean distance between the boundary components of the annulus
${\Cal V}^+\setminus E$. Then we have 
$$
\delta\ge |J|\min\{\varepsilon_1, \varepsilon\}
\ .
\slabel{delj}
$$
Moreover, by condition G$_2$,
$$
\diam ({\Cal V}^+\setminus E)\le \diam ({\Cal V})\le K|J|
\ .
\slabel{diams}
$$

Combining {\eqnodelj} and {\eqnodiams} with inequalities {\eqnoone}
and {\eqnotwo} of Lemma {4.3}, we get an upper bound for the
diameter of $E$ in the hyperbolic metric of ${\Cal V}^+$, namely
$$
\diam_{{\Cal V}^+}(E)\;\le\; \frac{C_0\pi
  K^2}{4\left[\min\{\varepsilon_1,\varepsilon\}\right]^2} 
\;=\;s
\ .
\slabel{upper}
$$

Now let $z\in E$ and suppose that $F(z)=\xi(z)$, {\it i.e.\/} that
$z\in {\Cal O}_\xi^+$ (other cases being treated in the same way). If
$w$ is any point in the non-empty set $\partial{\Cal
O}_\xi^+\cap\partial E$, we have
$$
d_{{\Cal V}^+}(z,{\Cal V}^+\setminus {\Cal O}_\xi^+)\le d_{{\Cal
V}^+}(z,w)\le \diam_{{\Cal V}^+}(E) 
\ .
\slabel{distz}
$$
Viewed as a (univalent) map from ${\Cal V}^+$ into itself, $\xi^{-1}$
is the composition of the hyperbolic isometry $\xi^{-1}:{\Cal V}^+\to
{\Cal O}_\xi^+$ with the hyperbolic contraction given by the inclusion
${\Cal O}_\xi^+\subseteq {\Cal V}^+$. Therefore, by Lemma {4.2} and
{\eqnodistz},
$$
\|(\xi^{-1})'(\xi(z))\|_Y\le C\left(\diam_{{\Cal V}^+}(E)\right)<1
\ ,
$$
which put together with {\eqnoupper} gives us
$$
\|F'(z)\|_Y=\|\xi'(z)\|_Y\ge\frac{1}{C(s)}>1
\ .
$$
This proves part ($a$) with $\varrho(\varepsilon ,K)=C(s)^{-1}-1>0$.

To prove part ($b$) we note, using parts ($c$) and ($d$)
of Lemma {4.5}, that if $z\in \dom (\Gamma )\cap V(J,\varepsilon)$ is such
that $F(z)\notin V(J,\varepsilon)$, then $z\notin
V(J,\tau(\varepsilon,K))$. In particular, $\Im{z}\ge
\tau(\varepsilon,K)|J|$, and therefore by ($a$) we have
$$
\|F'(z)\|_Y\ge
1+\varrho(\tau(\varepsilon,K),K)=\lambda(\varepsilon,K) 
\ .
$$
This completes the proof of both assertions for $\|\cdot\|_Y$. The proof for
$\|\cdot\|_{Y^*}$ is similar.\qed
\enddemo



We now turn to our second expansion result.

\proclaim{Proposition {4.8}} Let $\Gamma$ be a $K$-controlled holomorphic
pair with arbitrary rotation number. Then for each $\mu>1$ there exists
$n_*=n_*(\mu,\Gamma)>0$ with the following property. If $z$ is a point in the
domain of $\Gamma$ such that $F^j(z)\in \dom (\Gamma)$ for all $0\le j<n$ but 
$F^n(z)\notin \dom (\Gamma)$ for some $n>n_*$, then $\|(F^n)'(z)\|_Y\ge\mu$.
\endproclaim

\demo{Proof} 
Let $\Delta=[\xi(0),\eta(0)]$ be the short dynamical interval of
$\Gamma$, and let $f:\Delta\to \Delta$ be the associated circle map.
Also, let $K_0$ and $n_0=n_0(\Gamma)$ be given by Theorem
{4.1}. Then for all $n\ge n_0$ the $n$-th renormalization of $\Gamma$
is a holomorphic pair with the same co-domain ${\Cal V}$ that
restricts to a $K_0$-controlled holomorphic pair
$\Gamma_n$. We denote by ${\Cal U}_n=\dom (\Gamma_n)$ the domain and
${\Cal V}_n$ the co-domain of $\Gamma_n$ , by $F_n$ its shadow, by
$J_n$ its long dynamical interval and by $I_n$ its short dynamical
interval. We also consider the protection bands ${\Cal
W}_n=V(J_n,\varepsilon_0)$, where $\varepsilon_0$ is the constant in Lemma
{4.5}. Recall that $F_n({\Cal W}_n)\subseteq \dom(\Gamma_n)$.

Since $|J_n|\to 0$ as $n\to\infty$, and since $\diam {\Cal
V}_n\le K_0|J_n|$ by condition G$_2$, there exists a sequence $n_0<
n_1<n_2<\cdots <n_j<\cdots$ along which domains and co-domains are
nested in the sense that  
$$
{\Cal V}_{n_{j+1}}\subseteq {\Cal U}_{n_j}\cap {\Cal W}_{n_j}
\slabel{nest}
$$ 
holds true for all $j\ge 1$. Now choose an integer $k$ such that
$\lambda^{k-1}>\mu$, where $\lambda=\lambda(\varepsilon_0,K_0)$ is 
given by Proposition {4.7}~($b$).
\bigskip
$$
\psannotate{
\psboxto(12cm;0cm){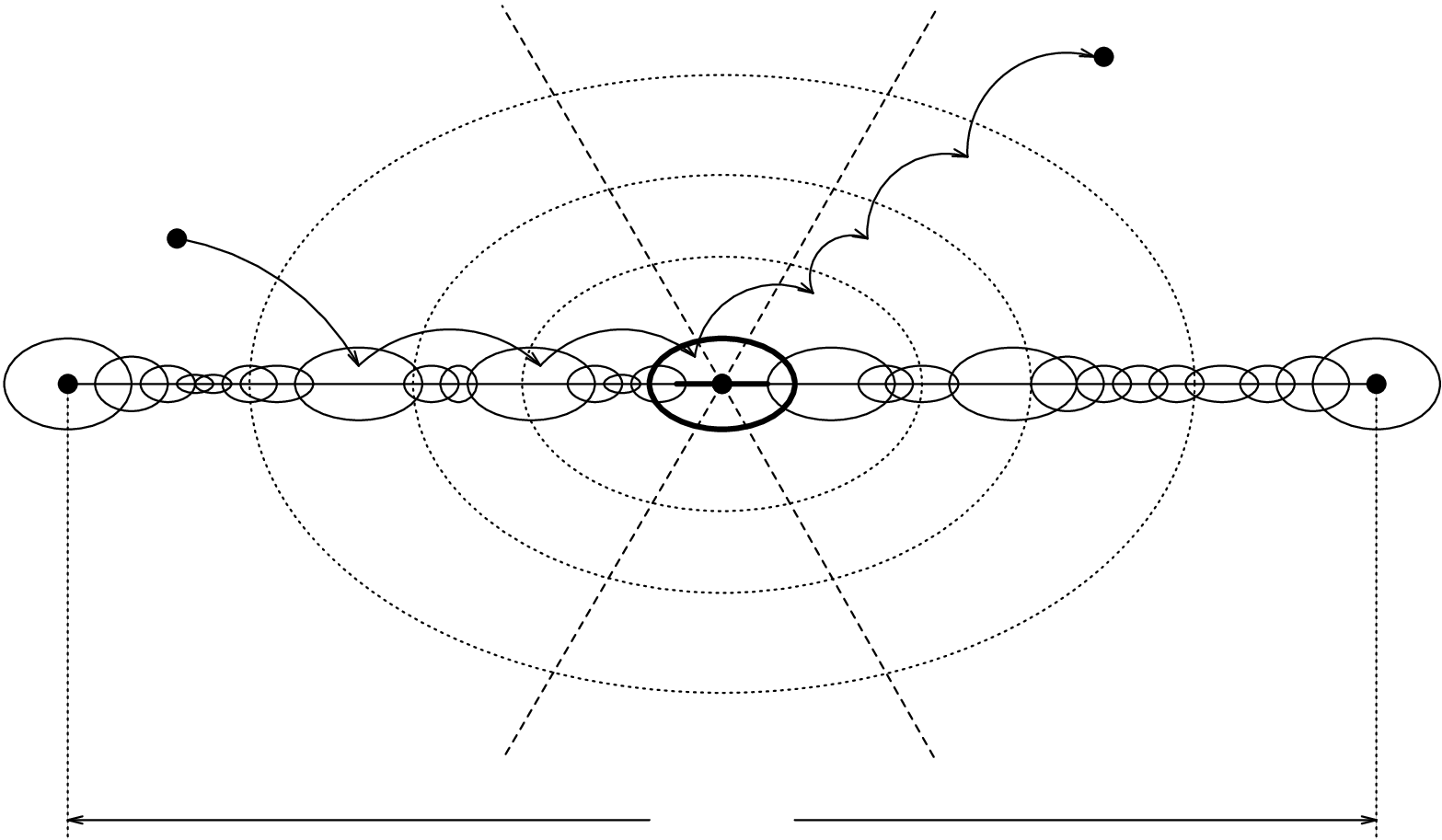}}
{\fillinggrid\at{13.7\pscm}{0.1\pscm}{$J$}
             \at{3.5\pscm}{10.9\pscm}{$w$}
             \at{13.2\pscm}{6.5\pscm}{${\Cal U}_{n_k}$}
             \at{20.1\pscm}{14\pscm}{$F_{n_j}^{s_j}(z_j)$}
             \at{14.8\pscm}{5.45\pscm}{${\Cal U}_{n_{k-1}}$}
             \at{18.8\pscm}{3.4\pscm}{${\Cal U}_{n_j}$}
}
$$
\smallskip
\centerline{\it Figure 3}
\bigskip
{\it Claim.\/} If $z\in {\Cal U}_{n_k}$ and $n>0$ are such that
$F^j(z)\in \dom (\Gamma)$ for $j=0,1,\ldots ,n-1$ but $F^n(z)\notin
\dom (\Gamma)$, then $\|(F^n)'(z)\|_Y>\mu$.
\smallskip
To prove this claim, we proceed as follows. Since $z\in {\Cal
  U}_{n_k}$ and the $F$-orbit of $z$ eventually exits the domain of
$\Gamma$, there exists $s_k\ge 1$ such that $F_{n_k}^i(z)\in {\Cal
  U}_{n_k}$ for all $0\le i<s_k$ but $z_{k-1}=F_{n_k}^{s_k}(z)\notin {\Cal
  U}_{n_k}$. By the nesting condition {\eqnonest}, we must have
$z_{k-1}\in {\Cal U}_{n_{k-1}}\cap {\Cal W}_{n_{k-1}}$. Proceeding inductively
  in this fashion until the $F$-orbit of $z$ exits ${\Cal U}_{n_1}$, we obtain
  positive integers $s_k,s_{k-1},\ldots ,s_1$ and points $z_k=z,
  z_{k-1},\ldots ,z_0$ satisfying
$$
z_{j-1}=F_{n_j}^{s_j}(z_j)\in {\Cal U}_{n_{j-1}}\cap {\Cal W}_{n_{j-1}}
\ ,
$$
and such that $F_{n_j}^i(z_j)\in {\Cal U}_{n_j}$ for all $0\le i<s_j$
while $z_{j-1}\notin {\Cal U}_{n_j}$, for all $1\le j\le k$. Note that since
the $F$-orbit of $z$ exits each ${\Cal U}_{n_j}$, it exits each protection
band ${\Cal W}_{n_j}$ {\it a fortiori\/}. Therefore, by
Proposition {4.7}~($b$), we have $\|(F_{n_j}^{s_j})'(z_j)\|_Y\ge \lambda$ for
$j=1,2,\ldots,k-1$. When the $F$-orbit of $z$ reaches $z_0$, there are 
still, say, $r$ iterates to go before it finally exits the domain of
$\Gamma$, and we see that 
$$
F^n(z)\;=\;F^r\circ F_{n_1}^{s_1}\circ F_{n_2}^{s_2}\circ\cdots \circ
F_{n_k}^{s_k}(z)
\ .
$$
Therefore by the chain rule we have
$$
\|(F^n)'(z)\|_Y\;=\;\|(F^r)'(z_0)\|_Y\,\prod_{j=1}^k
\|(F_{n_j}^{s_j})'(z_j)\|_Y\;>\; 
\lambda^{k-1}>\mu
\ ,
$$
which proves the claim.

The next step ({\it cf.\/} Figure 3) is to cover $J$ with pre-images 
of ${\Cal U}_{n_k}$ under $F$. More precisely, since ${\Cal U}_{n_k}
\supseteq I_{n_k}$ and  
$$
\Delta\;=\; \bigcup_{j=0}^{q_{n_k}}f^{-j}(I_{n_k})
\ ,
$$
it follows that $\Delta$ is contained in the open set
$$
{\Cal O}=\bigcup_{j=0}^{q_{n_k}}F^{-j}({\Cal U}_{n_k})
\ .
$$
Here we have used the fact that $f$ is simply the restriction of $F$
to the short dynamical interval of $\Gamma$. Now recall that
$J=\Delta\cup \eta^{-1}(\Delta)\cup \xi^{-1}(\Delta)\cup\cdots
\cup\xi^{-m}(\Delta)$, where $m$ is the {\it height} of $\Gamma$.
Hence $J$ is contained in the open set  
$$
{\Cal O}'\;=\; \bigcup_{j=0}^m F_*^{-j}({\Cal O})
\ ,
$$
where $F_*$ is the extended shadow of $\Gamma$ (we use $F_*$ instead
of $F$ so we can cover the endpoints of $J$ also).
This means that for some $\varepsilon>0$ small enough
$V(J,\varepsilon)\subseteq {\Cal O}'$. 
All points in the domain of $\Gamma$ that are not in $V(J,\varepsilon)$ lie at
a definite positive distance from the real axis, in other words,
$$
\varepsilon'=\inf\,\big\{|\Im{z}|:z\in\dom(\Gamma)\setminus
V(J,\varepsilon)\big\}\;>\;0 
\ .
$$
Hence, let $n_*$ be such that $(1+\varrho)^{n_*}>\mu$, where $\varrho=\varrho
(\varepsilon', K)$ is given by Proposition {4.7}~($a$). If $w\in \dom
(\Gamma)$ and $n>n_*$ are such that $F^j(w)\in \dom (\Gamma)$ for all $0\le
j<n$ but $F^n(w)\notin \dom (\Gamma)$, then {\it either} there exists $j$ such
that $F^j(w)\in V(J,\varepsilon)$, in which case there exists $i\ge j$
such that $z=F^i(w) \in {\Cal U}_{n_k}$ and therefore by the Claim above 
$$
\|(F^n)'(w)\|_Y\ge\|(F^{n-i})'(z)\|_Y\ge \mu
\ ,
$$
{\it or} $F^j(w)\notin
V(J,\varepsilon)$ for all $j<n$, in which case
$|\Im{F^j(w)}|\ge\varepsilon'|J|$, whence by Proposition {4.7}~($a$)
and the chain rule we get $\|(F^n)'(w)\|_Y\ge (1+\varrho)^n>\mu$ also.
This completes the proof. \qed
\enddemo

\proclaim{Proposition {4.9}} Let $\Gamma$ be a holomorphic pair with
arbitrary rotation number, and let ${\Cal K}_\Gamma$ be its limit
set. Then 
$$
{\Cal K}_\Gamma=\overline{\bigcup_{n\ge 0}  F^{-n}(J)}
\ .
$$
\endproclaim
\demo{Proof} Since the statement we want to prove is purely
topological, we can assume, using the pull-back theorem for
holomorphic pairs (Theorem {2.3}) and Theorem {4.1},
that $\Gamma$ satisfies all the hypotheses of Proposition {4.8}. 
For the same reasons, we also know that $\Gamma$ is
conjugate to a deep renormalization of a suitable element $f$ of the
Arnold family. More precisely, there exist domains ${\Cal V}^0, 
{\Cal O}_\xi^0, {\Cal O}_\eta^0, {\Cal O}_\nu^0\subseteq {\Bbb C}^*$,
symmetric about $\partial {\Bbb D}$ with respect to inversion, forming
a (generalized) bowtie with center at $1\in \partial {\Bbb D}$, and a
positive integer $s$ such that the maps
$$
\xi_0=f^{q_s}|{\Cal O}_\xi^0,\ 
\eta_0=f^{q_{s+1}}|{\Cal O}_\eta^0,\ 
\nu_0=f^{q_s+q_{s+1}}|{\Cal O}_\nu^0 \ 
$$
determine a (generalized) holomorphic pair $\Gamma^0$, and moreover
there exists a quasi-conformal homeomorphism $H:{\Cal V}\to {\Cal
  V}^0$ conjugating $\Gamma$ to $\Gamma^0$. Let $F_0$ be the shadow of
$\Gamma^0$, and let $J_0=H(J)\in \partial {\Bbb D}$ be the long
dynamical interval of $\Gamma^0$. The distinctive feature of
$\Gamma^0$ is the fact that, for each $w\in \dom (\Gamma^0)$ and each
$k>0$ for which $F_0^k(w)$ is defined, there exist non-negative
integers $m,n$ with $k\le m+n$ such that 
$$
F_0^k(w)\;=\;f^{mq+nQ}(w)
\ ,
\slabel{word}
$$
where $q=q_s$ and $Q=q_{s+1}$.
\smallskip
{\it Claim.\/} The set $\Lambda={\Cal K}_\Gamma\setminus\bigcup_{n\ge
  0}F^{-n}(J)$ has empty interior.
\smallskip
To prove this, we argue by contradiction. Suppose $D\subseteq \Lambda$
is a non-empty open disk, and let $D_0=H(D)$. Since points in $D$ can
be iterated by $F$ forever, $F_0^k(w)$ exists for all $k\ge 0$, for
each $w\in D_0$, and so by {\eqnoword} there exist infinitely
many $n>0$ such that $f^n(w)\in \dom (\Gamma_0)$. We know that the
Julia set of $f$ is ${\Bbb C}^*$, so by Montel's theorem there exist
$n_0>0$ and $w\in D_0$ such that $f^{n_0}(w)\in \partial {\Bbb D}$. But
then $f^n(w)\in \partial {\Bbb D}$ for all $n\ge n_0$, and so there
exists $k>0$ such that $F_0^k(w)\in \partial {\Bbb D}\cap\dom
(\Gamma_0)=J_0$. We have thus found a point $z=H^{-1}(w)$ in $D$ which 
belongs to $F^{-k}(J)$, a contradiction that proves the claim.

Now, let $x$ be an arbitrary point of $\Lambda$. We want to show that
for every $\delta>0$ there exists a point $p\in D(x,\delta)$ such that
$F^m(p)\in J$ for some $m>0$. Since $F^n(x)$ exists for all $n\ge 0$
and lies in the {\it interior} of $\dom (\Gamma)$, it follows from the
claim above that for all $\delta>0$ and all sufficiently large $k$
there exists $y\in D(x,\delta)$ such that $F^j(y)\in \dom (\Gamma)$
for $0\le j<k$ but $F^k(y)\notin \dom (\Gamma)$. By an easy
connectedness argument we can assume also that $F^k(y)\in \partial\dom
(\Gamma)$. We call such $y$ a {\it $(k,\delta)$-escaping point for $x$\/}. 

Next, observe that $C_1\;=\;\sup\,\{d_Y(z,F^{-1}(J)\cap Y):
z\in\partial\dom(\Gamma)\}<\infty$. Given $\delta>0$, let $\mu>1$ be such that
$$
\mu>\left(\frac{2C_1}{\delta}\right)^{e^{C_0C_1}}
\ ,
\slabel{muone}
$$
where $C_0$ is the universal constant of Lemma {4.4}. Let $y$ be a
$(n,\delta/2)$-escaping point for 
$x$ with $n>n_*$, where $n_*$ is given by Proposition {4.8}, and
let $w=F^n(y)\in\partial\dom (\Gamma)$. Take any point $w^*\in F^{-1}(J)$ such
that $d_Y(w,w^*)\le C_1$. Denote by $F^{-n}$ the inverse
branch of $F^n$ that maps $w$ back to $y$, and let $p=F^{-n}(w^*)$.  

To estimate the Euclidean distance between $x$ and $p$, we first
estimate the hyperbolic distance in $Y$ between $y$ and $p$. For this
purpose, let $\sigma$ be the hyperbolic geodesic in $Y$ joining $w$ to
$w^*$, and let $\rho_Y$ be the hyperbolic density of $Y$. The distance
$d_Y(y,p)$ is certainly not greater than the hyperbolic length
$L_Y(F^{-n}\sigma)$ of the arc $F^{-n}\sigma$ joining $y$ to $p$. On
the other hand, by Proposition {4.8} we have $\|(F^{-n})'(w)\|_Y\le
\mu^{-1}$, and applying Lemma {4.4} we get
$$
\|(F^{-n})'(z)\|_Y\le \frac{1}{\mu^{1/\alpha}}
\ ,
\slabel{mutwo}
$$
for all $z\in \sigma$, where $\alpha=\exp\{C_0d_Y(w,z)\}\le
e^{C_0C_1}$. Hence, using {\eqnomuone} and {\eqnomutwo}, we see
that 
$$
\eqalign{
L_Y(F^{-n}\sigma)&=\int_\sigma\|(F^{-n})'(z)\|_Y\,\rho_Y(z)|dz|\cr
&{}\cr
&\le \frac{1}{\mu^{1/\alpha}}\,\int_\sigma \rho_Y(z)|dz| =
\frac{1}{\mu^{1/\alpha}}\,d_Y(w,w^*)\cr 
&{}\cr
&\le C_1\mu^{-e^{-C_0C_1}}< \frac{\delta}{2}\ .\cr}
$$
This gives us $|y-p|\le d_Y(y,p)<\delta/2$, and therefore
$|x-p|<\delta$ as desired. \qed
\enddemo

The main consequence of Proposition {4.9} is one of the key points
we shall use in \S 6 to prove that the limit set of a
holomorphic pair is uniformly twisting. Let $z$ be a point in the limit
set ${\Cal K}_\Gamma$ of a holomorphic pair $\Gamma$. For each tangent
vector $v$ at $z$, we write $\ell_Y(v)$ for the length of $v$ measured
in the hyperbolic metric of $Y$, and put
$\ell_Y(v)=\infty$ if $z$ happens to lie in the real axis. For each
$k\ge 0$, we let $v_k=DF^k(z)\,v$.

\proclaim{Corollary {4.10}} Let $\Gamma$ be a holomorphic pair with
arbitrary rotation number, and let $z$ be a point in the limit set
${\Cal K}_\Gamma$. Then, for each tangent vector $v$ at $z$, we have
$$
\ell_Y(v_0)\le \ell_Y(v_1)\le \cdots \le\ell_Y(v_k)\;\to\;\infty
$$
as $k\to\infty$.
\endproclaim

\demo{Proof} The sequence is clearly non-decreasing, so we concentrate
on proving that it diverges to infinity. Assume $z\in {\Cal
  K}_\Gamma\cap Y$, otherwise there is nothing to prove, and fix
$k>0$. Let ${\Cal O}_k\subseteq \dom (\Gamma)$ be the connected
component of $F^{-k}(Y)$ that contains $z$. We know that $F^k$ maps
${\Cal O}_k$ univalently onto one of the two connected components of
$Y$. For definiteness, let us assume that ${\Cal O}_k\subseteq {\Cal
  V}^+$ and that $F^k({\Cal O}_k)={\Cal V}^+$. Thus, viewed as a map
from ${\Cal V}^+$ into itself, $F^{-k}$ is the composition of an
{\it isometry} between ${\Cal V}^+$ and ${\Cal O}_k$, with their
respective hyperbolic metrics, and a {\it contraction} given by the
inclusion ${\Cal O}_k\subseteq {\Cal V}^+$. In particular, 
$$
\|DF^{-k}(F^k(z))\|_Y\;=\;\frac{\rho_Y(z)}{\rho_{{\Cal O}_k}(z)} 
\ .
\slabel{ell}
$$
Now, if $D_k$ is the largest disk centered at $z$ and contained in
${\Cal O}_k$, then
$$
\rho_{{\Cal O}_k}(z)\;\asymp\;\rho_{D_k}(z)\;=\;\frac{2}{R_k}
\ ,
$$
where $R_k$ is the radius of $D_k$. We also saw in Lemma {4.6}
that $\rho_Y(z)\asymp 1/|\Im{z}|$. Combining these facts 
with {\eqnoell} and the chain rule, we get
$$
\ell_Y(v_k)\;=\;\|(F^k)'(z)\|_Y\,\ell_Y(v)\asymp \frac{2|\Im{z}|}{R_k}
\,\ell_Y(v) 
\ .
$$
Since by Proposition {4.9} we have $R_k\to 0$ as $k\to\infty$,
we are done.\qed
\enddemo

\heading
5. The critical point is deep
\endheading

Now we use the expansion results established in the previous section
to prove that the critical point of a holomorphic pair $\Gamma$ is a
$\delta$-deep point of the limit set ${\Cal K}_\Gamma$ of $\Gamma$, for some
$\delta>0$. By definition, this means that for every $r>0$ the largest
disk contained in $D(0,r)$ which does not intersect ${\Cal K}_\Gamma$ has
radius $\le r^{1+\delta}$. (Thus, linear blow-ups of ${\Cal K}_\Gamma$ around
$0$ fill-out the plane at an exponential rate in the Hausdorff
metric.)

\proclaim{Theorem {5.1}} Let $\Gamma$ be a holomorphic pair with
arbitrary rotation number and limit set ${\Cal K}_\Gamma$. Then there exists
$\delta>0$ such that the critical point of $\Gamma$ is a $\delta$-deep point
of ${\Cal K}_\Gamma$. 
\endproclaim
\demo{Proof} 
In the course of the proof, we write $C_0$ for the constant of Lemma
{4.4}, and $C_1, C_2,\ldots$ for constants that depend only on the
complex bounds.
Let $\Gamma_n$ be the $n$-th renormalization of
$\Gamma$. By Theorem {4.1} we may assume that $\Gamma_n$ is
$K_0$-controlled for all $n\ge n_0$. As in the proof of
Proposition {4.8}, we choose $n_0<n_1<\cdots n_i<\cdots$ so that ${\Cal
V}_{n_{i+1}}\subseteq {\Cal U}_{n_i}\cap {\Cal W}_{n_i}$, where the ${\Cal
W}_n$ are the protection bands defined in that Proposition. Since $\diam({\Cal
U}_n)\asymp \diam({\Cal V}_n)\asymp |J_n|$ and $|J_n|\to 0$ exponentially as 
$n\to\infty$, we can choose the $n_i$ so that $n_{i+1}-n_i$ is bounded
(by a constant depending only on the real bounds). 

We want to show that for every $z$ in the domain of $\Gamma$ such that
$z\notin {\Cal K}_\Gamma$ we have $\dist(z,{\Cal K}_\Gamma)\le
C_1|z|^{1+\delta}$ for some $\delta>0$. Let $k$ be the largest
such that $z\in {\Cal U}_{n_k}$. Then $\diam({\Cal U}_{n_k})\asymp
|z|$, and by the real bounds
$$
n_k\;\asymp\;\log{\left(\frac{1}{|J_{n_k}|}\right)}
\;\asymp\;\log{\left(\frac{1}{|z|}\right)} 
\ .
$$
But $n_k\le C_2k$, and we get $k\ge C_3\log{(1/|z|)}$. Now, since
$z\notin {\Cal K}_\Gamma$, its forward orbit eventually leaves ${\Cal
U}$, and so it leaves each ${\Cal U}_{n_i}$, and each ${\Cal W}_{n_i}$ as
well. Let $n$ be such that $F^n(z)\in {\Cal U}_{n_0}\cap {\Cal W}_{n_0}$ but
$F^{n+1}(z)\notin {\Cal W}_{n_0}$. Then, by the same argument used in the
proof of the Claim in Proposition {4.8}, we have
$$
\|DF^n(z)\|_Y\;\ge\;\lambda^k\;\ge\;
\left(\frac{1}{|z|}\right)^{C_3\log{\lambda}} 
\;=\;\frac{1}{|z|^{\beta}}
\ .
$$
Moreover, by Lemma {4.5} ($c$), the point $z'=F^n(z)$ satisfies
$\Im{z'}\ge\tau(\varepsilon_0,K_0)|J_{n_0}|$, and so we can find a point
$w'\in {\Cal K}_{\Gamma}\cap \partial {\Cal W}_{n_0}$ such that $d_Y(z',w')\le
C_4$. Let $\gamma'$ be the geodesic arc in $Y$ joining $z'$ to $w'$ (so that
$L_Y(\gamma')\le C_4$). Then the arc 
$\gamma=F^{-n}(\gamma')$ joins $z$ to $w=F^{-n}(w')\in {\Cal K}_\Gamma$. 
Using Lemma {4.4}, we see that $\|DF^{-n}(\zeta)\|_Y\le
\|DF^{-n}(z)\|_Y^\alpha$ for all $\zeta\in \gamma'$, where
$\alpha=e^{C_0C_4}$. Therefore $L_Y(\gamma)\le C_5|z|^{\delta}$, where
$\delta=\alpha\beta$. Since the hyperbolic density of $Y$ at points in ${\Cal
U}$ is comparable to the hyperbolic density of the upper half-plane ${\Bbb
H}$, it follows that $L_{{\Bbb H}}(\gamma)\le C_6|z|^\delta$. Thus, $w$
belongs to the hyperbolic disk $\Omega\subseteq {\Bbb H}$ of center $z$ and
radius $C_6|z|^\delta$. Since the point of $\partial\Omega$ farthest from $z$ 
in the Euclidean metric lies vertically above $z$, let us say $z+iL$,
we get
$$
\log{\left(1+\frac{L}{\Im{z}}\right)}\;\le\;C_6|z|^\delta
$$
But since $|z|$ is bounded, the left-hand side of this inequality is
$\ge C_7L/\Im{z}$, and we finally have
$$
|w-z|\;\le\;L\;\le\;\frac{C_6}{C_7}|z|^\delta(\Im{z})\;\le\;C_8|z|^{1+\delta}
\ .
$$
We deduce that the largest Euclidean disk centered at $z\notin {\Cal
K}_\Gamma$ which does not meet ${\Cal K}_\Gamma$ has radius $\le
C_8|z|^{1+\delta}$, and this proves the theorem. \qed
\enddemo

\heading
6. Small limit sets everywhere
\endheading

In this section we use Corollary {4.10} to prove an important
property of holomorphic pairs that will imply (\S 7) that the limit
set is {\it uniformly twisting} in the sense of McMullen {\refMcb}.
The precise statement is given in Theorem {6.3} below. We will
need the following property of critical circle maps.

\proclaim{Lemma {6.1}} Let $f:S^1\to S^1$ be a critical circle map
for which the real a priori bounds hold true. For each $n\ge 1$, let
$\Delta_n$ be the interval of endpoints $f^{q_n}(c)$ and $f^{q_{n-1}}(c)$
containing the critical point $c$. Then $\Delta_n, f^{-1}(\Delta_n), \ldots,
f^{-q_n+1}(\Delta_n)$ are pairwise disjoint and each $f^{-j}(\Delta_n)$ has
definite and bounded space on both sides inside the largest interval
containing it that does not meet any of the others. 
\endproclaim
\demo{Proof} The endpoints of $f^{-j}(\Delta_n)$ are $f^{q_n-j}(c)$ and
$f^{q_{n+1}-j}(c)$, and both belong to the dynamical partition ${\Cal
P}_n$ of level $n$ of $f$.\qed
\enddemo

We will also need an easy property of controlled holomorphic pairs. 

\proclaim{Lemma {6.2}} Let $\Gamma$ be a controlled holomorphic
pair. Then there exist a disk $D'\subseteq {\Cal O}_\nu\setminus {\Bbb R}$
and a disk $D''\subseteq {\Cal U}$ centered at the origin, with
$\diam(D')\asymp |J|\asymp \diam(D'')$ and $\dist(D',{\Bbb R})\asymp
|J|$, such that $\nu$ is univalent in $D'$ and $\nu(D')\supseteq
D''$.\qed 
\endproclaim

\demo{Definition} We say that a holomorphic pair $\Gamma$ is {\it
super-controlled} if $\Gamma$ is $K$-controlled and there
exist $R=R(K)>0$, domains $\widetilde{\Cal
O}_\xi,\widetilde{\Cal O}_\eta,\widetilde{\Cal O}_\nu$, and
complex-analytic extensions $\widetilde\xi: \widetilde{\Cal O}_\xi\to
{\Bbb C}$, $\widetilde\eta: \widetilde{\Cal O}_\eta\to {\Bbb C}$,
$\widetilde\nu: \widetilde{\Cal O}_\nu\to {\Bbb C}$ of $\xi$, $\eta$,
$\nu$ respectively, such that for $\gamma=\xi,\eta,\nu$ and for all
$z\in \partial {\Cal O}_\gamma\cap {\Cal K}_\Gamma$ we have
$D(z,R|J|)\subseteq \widetilde{\Cal O}_\gamma$ and the restriction of
$\widetilde\gamma$ to $D(z,R|J|)$ has distortion bounded by $K$. 
\enddemo

To go from control to super-control, one renormalization suffices. 
In other words, if $\Gamma$ is $K$-controlled then its first renormalization
${\Cal R}(\Gamma)$ is a super-controlled holomorphic pair.

\proclaim{Theorem {6.3}} Let $\Gamma$ be a super-controlled
holomorphic pair with rotation number of bounded type, with limit set ${\Cal
K}_\Gamma$ and shadow $F$. Then 
for each $z_0\in{\Cal K}_\Gamma$ and each $r>0$ there exists a pointed
domain $(U,y)$ with $|z_0-y|=O(r)$ and $\diam(U)\asymp r$, and there
exists $k>0$ such that $F^k$ maps $(U,y)$ onto a pointed domain
$(V,0)$ univalently with bounded distortion. In particular, $U$
contains a conformal copy of a high renormalization of $\Gamma$ whose
limit set has size commensurable with $r$. 
\endproclaim
\demo{Proof} The proof is divided into several steps.
\smallskip
\noindent {\bf{Step I.}\enspace}
First we show that the statement is true when $z_0\in
[\xi(0),\eta(0)]$. Let $f=F|[\xi(0),\eta(0)]$ and for each $n\ge 1$
let $\Delta_n$ be as in Lemma {6.1}. Each $\Delta_n$ is the short
dynamical interval of the $n$-th renormalization $\Gamma_n$ of
$\Gamma$. Let ${\Cal V}_n\supseteq \Delta_n$ be the co-domain of
$\Gamma_n$; by the complex bounds, we can assume that
$\diam({\Cal V}_n)\asymp |\Delta_n|$. For each $0<i\le q_n-1$, let
$V_{n,i}=f^{-i}({\Cal V}_n)$; note that by Koebe's distortion theorem we have
$\diam(V_{n,i})\asymp |f^{-i}(\Delta_n)|$. Now, choose $n$ so that the
interval of ${\Cal P}_n$ which contains $z_0$ has length $\asymp r$.
Then choose $k$ so that $V_{n,k}$ is closest to $z_0$. 
This value of $k$ and the pointed domain $(U,y)=(V_{n,k},f^{-k}(0))$
will do. 
\smallskip
\noindent {\bf{Step II.}\enspace}
Next, the statement is true for all $z_0\in (a,b)$, the long dynamical
interval of $\Gamma$. Here, let $p$ be the {\it height} of $\Gamma$
and consider the domains $V_{n,i}$ together with $\xi^{-1}(V_{n,i}),
\xi^{-2}(V_{n,i}),\ldots \xi^{-p}(V_{n,i})$ and $\eta^{-1}(V_{n,i})$.
Then apply the same argument in step I.
\smallskip
\noindent {\bf{Step III.}\enspace}
Now we prove the more difficult case when $z_0$ is not on the real
axis. Here is where we use Corollary {4.10}. This case breaks down
into further sub-cases. We start with a vector $v_0$ at $z_0$ and
Euclidean norm $|v_0|\asymp r$, and we iterate: $z_k=F^k(z)$,
$v_k=DF^k(z)v$. Note that $|F'(z_k)|=|v_{k+1}|/|v_k|$ for all $k$. We
also know from Corollary {4.10} that $\ell_Y(v_k)\to \infty$
monotonically as $k\to \infty$. There are two possibilities.
\smallskip
{\it Step IIIa.\enspace} There exists $k$ such that $\ell_Y(v_k)<< 1$
but $\ell_Y(v_{k+1})>> 1$. More precisely, $\ell_Y(v_k)<\varepsilon$
while $\ell_Y(v_{k+1})>1/\varepsilon$, for some constant $\varepsilon$
that will be determined in the course of the argument.
To be definite, assume $z_k$ lies in the domain of $\xi$, so
$z_{k+1}=\xi(z_k)$. Since at all points $z$ of the domain of $\Gamma$
the hyperbolic density of $Y$ is comparable to $1/\Im{z}$, we have
$$
\frac{|v_k|}{\Im{z_k}}\;\asymp\;\ell_Y(v_k)\;<\;\varepsilon
\ .
$$
Hence $\Im{z_k}\ge C|v_k|/\varepsilon >> |v_k|$. This means that $\xi$
is univalent in a disk $D(z_k,R)$ of radius $R\asymp\Im{z_k}$ (provided
it is well-defined there -- it may happen that $z_k$ is too close to
the boundary of the domain of $\xi$, in which case we must remember
that $\xi$ extends holomorphically to a definite neighborhood of such
boundary because the holomorphic pair $\Gamma$ is super-controlled). By
Koebe's one-quarter theorem, $\xi(D(z_k,R))$ contains the disk
$D(z_{k+1},R')$, where $R'=R|v_{k+1}|/4|v_k|>>|v_{k+1}|$. But
$\ell_Y(v_{k+1})\asymp |v_{k+1}|/\Im{z_{k+1}}$ is large
($>1/\varepsilon$), so $|v_{k+1}|>> \Im{z_{k+1}}$. Let
$\zeta=\Re{z_{k+1}}$ be the point on the real axis closest to
$z_{k+1}$. By step II, there exists a pointed domain $(U'',y'')$, which
is mapped univalently by some iterate of $F$ to a pointed domain $(V,0)$
around the critical point, such that $\diam(U'')\asymp |v_{k+1}|$ and
and $|\zeta-y''|=O(|v_{k+1}|)$, so $|z_{k+1}-y''|=O(|v_{k+1}|)$ also.
Now, if $\varepsilon$ is chosen small enough, we have $U''\subseteq
D(z_{k+1}, R'/2)$. Take $U'=\xi^{-1}(U'')\subseteq D(z_k,R)$ and
$y'-\xi^{-1}(y'')\in U'$. Again by Koebe's distortion theorem we have
$\diam(U')\asymp |v_k|$ and $|z_k-y'|=O(|v_k|)$. Since
$D(z_k,R)\subseteq Y$, the inverse branch $F^{-k}$ mapping $z_k$ back
to $z_0$ is well-defined and univalent in $D(z_k,R)$. Thus, take
$U=F^{-k}(U')$ and $y=F^{-k}(y')$. Once again by Koebe's distortion
theorem, $\diam(U)\asymp |v_0|\asymp r$, and
$|y-z_0|=O(|v_0|)=O(|r|)$, so we are done in this case.
\smallskip
{\it Step IIIb.\enspace} There exists $k$ such that $\ell_Y(v_k)\asymp
1$. As in McMullen's original argument, this case is the most delicate. 
Here, we have $|v_k|\asymp \Im{z_k}$. As before, let $\zeta=\Re{z_k}$
be the point on the real axis closest to $z_k$. Using step II, we find
some $V_{n,i}$ and $y_{n,i}\in V_{n,i}$ such that
$\diam(V_{n,i})\asymp \dist(\zeta, y_{n,i})\asymp \Im{z_k}$ (note that it is
here that we use the bounded type assumption) and such
that $F^i$ maps $(V_{n,i},y_{n,i})$ univalently onto $({\Cal V}_n,
0)$. Applying Lemma {6.2} to the controlled holomorphic pair
$\Gamma_n$, we get a disk $D'$. Let $U'=F^{-i}(D')\subseteq V_{n,i}$.
Using Koebe's distortion theorem, we see that $\diam(U')\asymp
\dist(U',{\Bbb R}) \asymp \diam(V_{n,i})$. It follows at once that if
$\gamma$ is the geodesic arc in $Y$ joining $z_k$ to $U'$, then
$\diam(\gamma \cup U')\asymp \Im{z_k}$. Hence we can find $r_k\asymp
\Im{z_k}$ and points $\zeta_1=z_k, \zeta_2,\ldots, \zeta_s\in
\gamma\cup U'$ such that 
$$
\gamma\cup U'\;\subseteq\; \Omega\;=\;\bigcup_{j=1}^s D(\zeta_j, r_k/2)
\ ,
$$
and such that $D(\zeta_j,r_k)\subseteq {\Cal U}\setminus {\Bbb R}$ for
all $j$. But now the inverse branch $\Phi=F^{-k}$ mapping $z_k$ back to
$z_0$ is well-defined over the union of the disks $D(\zeta_j,r_k)$.
Applying Koebe's distortion to each of these disks we see that for all
$z\in\Omega$ we have $|\Phi'(z)|\asymp |\Phi'(z_k)|\asymp
|v_0|/|v_k|$. But then we see that $U=\Phi(U')$ satisfies
$\diam(U)\asymp |v_0|=r$.
\qed
\enddemo

\heading
7. Proof of Theorem {1.1}
\endheading

Now, recall that a map $\phi:{\Bbb C}\to {\Bbb C}$ such that $\phi(0)=0$ is
said to be {\it $C^{1+\beta}$-conformal} at zero (for some $\beta>0$) if the
complex derivative $\phi'(0)$ exists and we have 
$\phi(z)=\phi'(0)z+O(|z|^{1+\beta})$ for all $z$ near zero. 

\proclaim{Theorem {7.1}} Let $\Gamma$ be a holomorphic pair with
rotation number of bounded type and let $H:\Gamma\simeq \Gamma'$ be a
quasiconformal conjugacy between $\Gamma$ and another holomorphic pair
$\Gamma'$. Then $H$ is $C^{1+\varepsilon}$-conformal at $0$, where
$\varepsilon>0$ is universal.
\endproclaim

The proof of this theorem follows McMullen's strategy in Ch.~9 of
{\refMcb}. We only need to show that our case falls within his
framework. The definition of holomorphic dynamical system adopted by
McMullen is sufficiently broad to include our case. Indeed, we can
define the {\it full dynamics} of a holomorphic pair $\Gamma$ as
follows. 

The relevant holomorphic relations are defined using the maps
$\xi:{\Cal O}_\xi\to {\Bbb C}$, $\eta:{\Cal O}_\eta\to {\Bbb C}$,
$\nu:{\Cal O}_\nu\to {\Bbb C}$ in their open domains.
Let $\Omega$ be the set of all {\it finite} words in the alphabet
$\{\xi,\eta,\nu\}$. Given $\omega\in \Omega$, let $\dom(\omega)$ be
the set of all $z\in \dom(\Gamma)$ such that $\omega(z)$ exists. Then
each $\dom(\omega)\subseteq {\Bbb C}$ is an open set. If
$\dom(\omega)$ is non-empty, we say that $\omega$ is admissible for
$\Gamma$, and the set of all such words is denoted by $\Omega_\Gamma$.
Now, for each pair $(\omega_1,\omega_2)\in \Omega_\Gamma\times
\Omega_\Gamma$, let
$$
G(\omega_1,\omega_2)\;=\;\big\{(z_1,z_2)\in\dom(\omega_1)\times
\dom(\omega_2): \omega_1(z_1)=\omega_2(z_2)\big\}\;
\subseteq {\widehat{\Bbb C}}\times {\widehat{\Bbb C}}
\ .
$$
Each $G(\omega_1,\omega_2)$ is an analytic hypersurface in
${\widehat{\Bbb C}}\times {\widehat{\Bbb C}}$. Now consider the set
${\Cal F}(\Gamma)$ of all holomorphic maps $g:U\to {\Bbb C}$ for which
there exists $(\omega_1,\omega_2)\in \Omega_\Gamma\times
\Omega_\Gamma$ such that ${\roman{gr}}(g)=\{(z,g(z)): z\in
U\}\subseteq G(\omega_1,\omega_2)$.  
This set is the {\it full dynamics} of $\Gamma$.

\demo{Proof of Theorem {7.1}} 
Using Theorem {6.3}, we prove that ${\Cal F}(\Gamma)$ is uniformly
twisting by the same argument as in the proof Theorem~9.18 in
{\refMcb}. Since by Theorem {5.1} the critical point of $\Gamma$ is
deep, we can now use Theorem~9.15 in {\refMcb}.\qed
\enddemo

\proclaim{Corollary {7.2}} Let $f$ and $g$ be real-analytic critical
circle maps with the same rotation number of bounded type, and let $h:f\simeq
g$ be the conjugacy that maps the critical point of $f$ to the critical point
of $g$. Then $h$ has a quasiconformal extension which is $C^{1+\varepsilon}$
conformal at the critical point. In particular, the successive
renormalizations of $f$ and $g$ converge together exponentially fast in the
$C^0$ sense.
\endproclaim
\demo{Proof} We know from Theorem {2.1} that $h$ is quasisymmetric. By
Theorem {3.1}, if $n$ is sufficiently large then the $n$-th
renormalizations ${\Cal R}^nf$ and ${\Cal R}^ng$ extend to holomorphic
pairs $\Gamma_n(f)$ and $\Gamma_n(g)$, respectively. By the pull-back 
theorem for holomorphic pairs (Theorem {2.3}), the restriction of
$h$ to the long dynamical interval of $\Gamma_n(f)$ extends to a 
quasiconformal conjugacy between $\Gamma_n(f)$ and $\Gamma_n(g)$. But
then Theorem {7.1} tells us that $H$ is $C^{1+\varepsilon}$
conformal at the critical point. The fact that ${\Cal R}^nf$ and
${\Cal R}^ng$ converge together at an exponential rate in the $C^0$
topology follows easily from this and Theorem {2.1}.\qed
\enddemo

The proof of Theorem {1.1} now follows from this corollary and the
fact, proved in {\refdFdM}, that exponential convergence of
renormalizations in the $C^0$-topology implies $C^{1+\alpha}$
conjugacy for some $0<\alpha<1$ depending on the rate of convergence, which is
universal. See also Th.~9.4 in {\refMS}, where a similar
result is proved in the context of unimodal maps.


\Refs
\widestnumber\key{MM8}

\ref \key\refnoAh
\by \Ah
\book Conformal invariants
\publ McGraw-Hill
\yr 1973
\endref

\ref \key\refnodF
\by \dF
\paper Asymptotic rigidity of scaling ratios for critical circle
mappings 
\publ to appear
\jour \ETDS
\endref

\ref \key\refnodFdM
\by \dF\ and \Me
\paper Rigidity of critical circle mappings I
\publ IMS Stony Brook preprint 97/16 (1997)
\endref


\ref \key\refnoGrS
\by \Gk\  and \Sw
\paper Critical circle maps near bifurcation
\yr 1996
\vol 176
\pages 227--260
\jour \CMP
\endref

\ref \key\refnoHea
\by \He
\paper Sur la conjugaison differentiable des diff\'eomorphismes du
cercle a des rotations 
\yr 1979
\vol 49
\pages 5--234
\jour \IHES
\endref

\ref \key\refnoHeb
\bysame
\paper Conjugaison quasi-sim\'etrique des hom\'eomorphismes du
cercle a des rotations 
\publ Manuscript
\yr 1988
\endref


\ref \key\refnoLYam
\by \Ly \ and \Yam
\paper Dynamics of quadratic polynomials: complex bounds for real maps
\publ MSRI Preprint 034-95
\yr 1995
\endref


\ref \key\refnoMca
\by \Mc
\paper Complex dynamics and renormalization
\jour Annals of Math Studies
\vol 135
\publ Princeton University Press
\yr 1994
\publaddr Princeton
\endref

\ref \key\refnoMcb
\bysame
\paper Renormalization and 3-manifolds which fiber over the circle
\jour Annals of Math Studies
\vol 142
\publ Princeton University Press
\yr 1996
\publaddr Princeton
\endref

\ref \key\refnoMcc
\bysame
\paper Self-similarity and of Siegel disks and Hausdorff dimension of Julia
sets 
\publ to appear
\jour \AcM
\endref


\ref \key\refnoMS
\by \Me \ and \vS
\book One dimensional dynamics
\publ Springer-Verlag
\yr 1993
\publaddr Berlin and New York
\endref

\ref\key\refnoSu
\by \Su
\paper Bounds, quadratic differentials and renormalization
conjectures 
\inbook Mathematics into the Twenty-First Century 
\bookinfo Amer. Math. Soc. Centennial Publication, vol.2 
\publ Amer. Math. Soc. \publaddr Providence, RI 
\yr 1991
\endref

\ref \key\refnoSwa
\by \Sw
\paper Rational rotation numbers for maps of the circle
\yr 1988
\vol 119
\pages 109--128
\jour \CMP
\endref

\ref \key\refnoYam
\by \Yam
\paper Complex bounds for critical circle maps
\publ IMS Stony Brook preprint 95/12
\yr 1995
\endref

\ref \key\refnoYo
\by \Yo
\paper Conjugaison analytique des diff\'eomorphismes du cercle
\publ Manuscript
\yr 1989
\endref




\endRefs
\enddocument